\documentclass[11pt]{article}
\usepackage{amssymb,amsmath,amsthm}
\usepackage{amscd}
\usepackage{mathrsfs}
\newcommand{\R}{\mathbb{R}}

\newcommand{\dd}{\,{\rm d}}

\newcommand{\e}{\mathbf{e}}

\numberwithin{equation}{section}
\newtheorem{thm}{Theorem}[section]
\newtheorem{df}[thm]{Definition}

\newtheorem{rem}[thm]{Remark}

\usepackage{verbatim}

\newcommand{\curl}{\text{curl}}

\usepackage{color}
\usepackage{enumitem}

\definecolor{Green}{rgb}{0.010,0.7,0.02}

\newcommand{\RR}{\mathbb{R}}
\newcommand{\pa}{\partial}

\begin{document}

\title{The Inviscid Limit and Boundary Layers for Navier-Stokes Flows}

\date{}

\author{
{\bf Yasunori Maekawa}\\
Department of Mathematics\\
Graduate School of Science\\
Kyoto University\\
Kitashirakawa Oiwake-cho, Sakyo-ku\\
Kyoto 606-8502, Japan\\
{\tt maekawa@math.kyoto-u.ac.jp}
\and
{\bf Anna Mazzucato}\\
Mathematics Department\\
The Pennsylvania State University\\
University Park, PA, 16802 USA\\
{\tt alm24@psu.edu}}

\maketitle

\begin{abstract}
The validity of the vanishing viscosity limit, that is, whether solutions of the Navier-Stokes equations modeling viscous incompressible flows converge to solutions of the Euler equations modeling inviscid incompressible flows as viscosity approaches zero, is one of the most fundamental issues in mathematical fluid mechanics. The problem is classified into two categories:
the case when the physical boundary is absent, and the case when the physical boundary is present and the effect of the boundary layer becomes significant.
The aim of this article is to review recent progress on the mathematical
analysis of this problem in  each category.\\
To appear in {\em Handbook of Mathematical Analysis in Mechanics of Viscous
Fluids},  Y. Giga and A. Novotn\'y Ed., Springer. The final publication is
available at www.springerlink.com.
\end{abstract}

\section{Introduction}

Determining the behavior of viscous flows at small viscosity is one of the most fundamental problems in fluid mechanics.
The importance of the effect of viscosity, representing tangential friction forces in fluids, is classical and has been recognized for a long time. A well-known example is the resolution of D'Alembert's paradox concerning the drag experienced by a body moving through a fluid, which is caused by neglecting the effect of viscosity in the theory of ideal fluids. For real flows like water and air, however, the kinematic viscosity is a very small quantity in many situation. Physically, the effect of viscosity is measured by a non-dimensional quantity, called the Reynolds number  $Re := \frac{UL}{\nu}$, where $U$ and $L$ are the characteristic velocity and length scale in the flow, respectively, and $\nu$ is a kinematic viscosity of the fluid. Therefore, when $U$ and $L$ remain in a fixed range, the limit of vanishing viscosity is directly related to the behavior of high-Reynolds number flows.
Hence, the theoretical treatment of the inviscid limit has great importance in
applications and has been pursued extensively in various settings. The aim of
this chapter is to give an overview of  recent progress in the inviscid limit
problem for incompressible Newtonian fluids,  although some results are also
available in the important cases of compressible fluids or non-Newtonian fluids.

The governing equations for incompressible homogeneous Newtonian fluids are the Navier-Stokes equations
\begin{equation} \label{e.NS}
\partial_t u^\nu + u^\nu \cdot \nabla u^\nu + \nabla p^\nu = \nu \Delta u^\nu +f^\nu \,, \qquad {\rm div}\, u^\nu=0\,,  \qquad \qquad u^\nu|_{t=0} = u_0^\nu\,.
\end{equation}
Here $u^\nu =(u^\nu_1,\cdots,u^\nu_n)$, $n=2,3$, and $p^\nu$ are the unknown velocity field and unknown pressure field, respectively, at time $t$ and position $x=(x_1,\cdots,x_n)\in \Omega$. In what follows, the domain $\Omega\subset \mathbb{R}^n$ will be either the whole space $\R^n$, or a domain with smooth boundary. Although many results can be stated for arbitrary dimension $n$, the physically relevant dimensions are $n=2,3$.
The external force $f^\nu=(f^\nu_1,\cdots,f^\nu_n)$ is a given vector field,
which will be typically taken as zero for simplicity, and
$u^\nu_0=(u^\nu_{0,1}\cdots,u^\nu_{0,n})$ is a given initial velocity field.
Standard notation will be used throughout for derivatives: $\partial_t =
\frac{\partial}{\partial t}$, $\partial_j = \frac{\partial}{\partial x_j}$,
$\nabla =(\partial_1, \cdots,\partial_n)$, $\Delta = \sum_{j=1}^n
\partial_j^2$, $u^\nu\cdot \nabla  =\sum_{j=1}^n u^\nu_j\partial_j$, and ${\rm
div}\, f = \sum_{j=1}^n \partial_j f_j$. The symbol $\nu$ represents the
kinematic viscosity of the fluid, which is taken as a small positive constant.
The equations \eqref{e.NS} are closed by imposing a suitable boundary condition,
which will be specified in different context (see Section \ref{sec.phys.b}). The
vorticity field is a fundamental physical quantity in fluid mechanics,
especially for incompressible flows, and it is defined as the curl of the
velocity field. Denoting $\omega^\nu={\rm curl}\, u^\nu$, one has:
\begin{align}
  \omega^\nu =\partial_1u_2^\nu -\partial_2 u^\nu_1 \quad (n=2)\,, \qquad \qquad
  \omega^\nu = \nabla \times u^\nu, \quad (n=3)
\end{align}
Above, we have identified  the vector $\omega^\nu=\varpi {\bf k}$, where ${\bf
k}=(0,0,1)$, with the scalar $\varpi$, and called the latter also $\omega^\nu$
with abuse of notation. For the reader's sake, we recall the
vorticity equations:
\begin{equation} \label{eq.vorticityEq}
   \pa_t \omega^\nu +u^\nu \cdot \nabla \omega^\nu =\nu \Delta \omega^\nu
    - \omega^\nu \cdot \nabla u^\nu. \qquad u^\nu = K_\Omega[ \omega^\nu],
\end{equation}
  where $K_\Omega$ stands for the Biot-Savart kernel in the domain $\Omega$. The last term on the right is called the {\em vorticity stretching term} and it is absent in two space dimensions. This term depends quadratically in $\omega^\nu$ and its presence precludes establishing long-time existence of Euler solutions in three space dimensions.

By formally taking the limit $\nu\rightarrow 0$ in \eqref{e.NS}, the Navier-Stokes equations are reduced to the Euler equations for incompressible flows
\begin{equation} \label{e.EE}
\partial_t u + u\cdot \nabla u + \nabla p = f\,, \qquad {\rm div}\, u=0\,, \qquad \qquad u|_{t=0} = u_0\,.
\end{equation}
The velocity field $u$ of the Euler flows will be often written as $u^0$ in this chapter.
Broadly speaking, a central theme of the inviscid limit problem is to understand when and in which sense the convergence of the Navier-Stokes flow $u^\nu$ to the Euler flow $u$ is rigorously justified.
Mathematically, the inviscid limit is a singular perturbation problem since the
highest order term $\nu\Delta u^\nu$ is formally dropped from the equations of
motion in the limit. For such a problem, in many cases the main issue is
establishing enough {\em a priori} regularity for the solutions  such that
convergence is guaranteed via a suitable compactness argument.  This singular
perturbation for the Navier-Stokes equations provides a challenging mathematical
problem, because of the fact that the non-linearity contains  derivatives of the
unknown solutions, and because of the non-locality arising from the pressure
term.
Indeed, even when the flow is two-dimensional and the fluid domain is the whole plane $\mathbb{R}^2$, the study of the inviscid limit problem becomes highly nontrivial if given data, such as initial data, possess little regularity.
It should be emphasized here that working  with non-smooth data has an important motivation not only mathematically but also in applications, for typical structures of concentrated vorticities observed in turbulent flows, such as the vortex sheets, vortex filaments, vortex patches, are naturally modeled as flows with certain singularities.
The problem for singular data but under the absence of physical boundaries will
be discussed in Section \ref{sec.phys.nb} of this chapter.

The inviscid limit problem becomes physically more important and challenging in
the presence of a nontrivial boundary, where the viscosity effects are found to
play a central role, in general, no matter how small the viscosity itself is,
depending on the geometry and on the boundary conditions for the flow.
Recent developments in the mathematical theory for this case will be reviewed in
Section \ref{sec.phys.b}.
The main obstruction in analyzing the inviscid limit in this context arises from
the complicated structure of the flow close to the boundary. Indeed, due to the
discrepancy between the boundary conditions in the Navier-Stokes equations and
in the Euler equations,
a boundary layer forms near the boundary where the effects of viscosity cannot
be neglected even at very low viscosity. The size and the stability property of
the boundary layer crucially depend on the type of prescribed boundary
conditions and also on the symmetry of the fluid domain and of the flows, which
are directly connected to the possibility of the resolution of the inviscid
limit problem. These are discussed in details in Sections \ref{s.slip} -
\ref{s.nonchar}. The concept of a viscous  boundary layer was first introduced
by Prandtl  \cite{Pra1904} in 1904 under the no-slip boundary condition, that
is, assuming  the flow adheres to the boundary. Since then, the theory of
boundary layers has had a strong impact in fluid mechanics
and has also initiated a fundamental approach in asymptotic analysis for
singular perturbation problems in differential equations. The reader is referred
to \cite{Schlichting} for various aspects of the boundary layer theory in fluid
mechanics. The basic idea is that the fluid region can be divided into two
regions:
a thin layer close to the boundary (the so-called boundary layer) where
the effect of viscosity is significant,
and the region outside this layer where the viscosity can be neglected and thus
the fluid behaves like an inviscid flow. As found  by Prandtl, the thickness of
the boundary layer is formally estimated as $\mathcal{O}(\sqrt{\nu})$,  at least
for no-slip boundary conditions, which is a natural scale given the parabolic
nature of the Navier-Stokes equations.
In the case of no-slip boundary conditions, the fundamental equations
describing the boundary layer are the Prandtl equations, which will be the
focus of Section \ref{s.Prandtl}.
Interestingly, the problem becomes most difficult in the case that the no-slip
boundary condition ,  the most classical and physically justified type
of boundary condition, is imposed,
 due to the fact that large gradients of velocity can then form at the boundary,
which may propagate in the bulk, giving rise to a strong instability mechanism
for the layer at high frequencies.
The rigorous description for the inviscid limit behavior of  Navier-Stokes
flows is still largely open in many relevant situations, and several recent
works in the literature have tackled this important problem.
In this chapter, standard notations in mathematical fluid mechanics will be used
for spaces of functions and vector fields. For example, the spaces
$C_{0,\sigma}^\infty (\Omega)$ and $L^p_\sigma (\Omega)$ are defined as
\begin{align*}
C_{0,\sigma}^\infty (\Omega) & = \{ f\in C_0^\infty (\Omega)^n ~|~ {\rm div}\, f=0~~{\rm in}~\Omega\}\,,\\
L^p_\sigma (\Omega) & = \overline{C_{0,\sigma}^\infty (\Omega)}^{\|f\|_{L^p(\Omega)}}\,, \qquad 1<p<\infty\,.
\end{align*}
The Bessel potential spaces $H^s(\mathbb{R}^n)$, the Sobolev spaces
$W^{s,p}(\Omega)$, $s\geq 0$, $1\leq p\leq \infty$, and the space of Lipschitz
continuous functions ${\rm Lip}\, (\mathbb{R}^n)$ are also defined as usual. The
$n$-product space $X^n$ will be often written as $X$ for simplicity.

\section{Inviscid limit problem without physical boundary}\label{sec.phys.nb}

This section is devoted to the analysis of the inviscid limit problem for the Navier-Stokes equations when the fluid domain has no physical boundary. In this case the effect of the boundary layer is absent at least from the physical boundary, and the problem is more tractable and has been analyzed in various functional settings.
Typically the Navier-Stokes flow is expected to converge to the Euler flow in the inviscid limit.
Then the main interest here is the class (regularity) of solutions for which this convergence is verified and its rate of convergence in a suitable topology.
To simplify the presentation, only the Cauchy problem in $\mathbb{R}^n$ will be the focus in this section, and therefore the external force in \eqref{e.NS} or \eqref{e.EE} will be  taken as zero.
To give an overview of known results, it will be convenient to classify the solutions depending on their regularity as follows.
\vspace{0.1cm}

\noindent
{\bf (I) Regular solutions}

\vspace{0.1cm}

\noindent
{\bf (II) Singular solutions}:\
(II-1) Bounded vorticity $\cdot$ Vortex patch ~~ (II-2) Vortex sheet $\cdot$ Vortex filament $\cdot$ Point vortices

\vspace{0.1cm}

\noindent
Since the external force is assumed to be zero, the above classification is essentially for the initial data.
Then, the typical case of ``(I) Regular solutions'' is that the initial data $u_0^\nu$ and $u_0$ belong to the Sobolev space $H^s(\mathbb{R}^n)$ with $s>\frac{n}{2}+1$, which is embedded in the space $C^1(\R^n)$.
On the other hand, the category ``(II) Singular solutions'' corresponds to the case when the initial data are less regular than the $C^1$ class.

\vspace{0.3cm}

\noindent
{\bf (I) Regular solutions.} The verification of the inviscid limit in $\mathbb{R}^n$ for smooth initial data (e.g., better than the $C^1$ class) is classical. Indeed, it is proved in  \cite{Go, Mc} for $\mathbb{R}^2$ and in \cite{Swann1971,Kato1972} for $\mathbb{R}^3$, and in \cite{EbMa} for  a compact manifold without boundary of any dimension.
The following result is given in \cite[Theorem 2.1]{Masmoudi2007}, and provides the relation between the regularity of solutions and the rate of convergence. For simplicity, the result is stated here only for the case $u^\nu_0=u_0$ in \eqref{e.NS} and \eqref{e.EE}, though the case $u^\nu_0\ne u_0$ is also given in \cite{Masmoudi2007}.

\begin{thm}\label{thm.regular} Let $u_0^\nu=u_0\in H^s (\mathbb{R}^n)$ with some $s>\frac{n}{2}+1$. Let $u\in C_{loc} ([0,T^*); H^s (\mathbb{R}^n))$ be the (unique) solution to the Euler equations, where $T^*>0$ is the time of existence of the solution. Then, for all $T \in (0,T^*)$ there exists $\nu_0>0$, such that for all $\nu\in [0,\nu_0]$ there exists a unique solution $u^\nu\in C([0,T]; H^s(\mathbb{R}^n))$ to the Navier-Stokes equations. Moreover, it follows that
\begin{align*}
\lim_{\nu\rightarrow 0} \| u^\nu  - u  \|_{L^\infty(0,T; H^s)} =0\,, \qquad \| u^\nu - u \|_{L^\infty(0,T; H^{s'})} \leq C(\nu T)^\frac{s-s'}{2}\,,
\end{align*}
for $s-2\leq s'\leq s$. Here, $C$ depends only on $u$ and $T$.
\end{thm}

Since the time $T^*$ in Theorem \ref{thm.regular} is the time of existence for the Euler flow, the inviscid limit  holds on any time interval in two space dimensions.
The  estimate in $H^{s'}(\mathbb{R}^n)$ for the case $s'=s-2$ is obtained by a standard energy method, and the case $s-2<s'<s$ is derived from interpolation. The convergence in $H^s(\mathbb{R}^n)$ is more delicate, and a regularization argument for the initial data is needed in the proof.

\vspace{0.3cm}

\noindent
{\bf (II) Singular solutions.} (II-1) Bounded vorticity $\cdot$ Vortex patch: The Euler equations are uniquely solvable, at least  locally in time, when the initial vorticity is bounded \cite{Y1963} or nearly bounded \cite{Yu2, Vi}. Hence it is expected that the Navier-Stokes flow converges to the Euler flow in the vanishing viscosity limit
also for such class of initial data. This class includes some important solutions, called vortex patches, which are typically vorticity fields defined as characteristic functions of bounded domains with smooth boundary.
Theorem \ref{thm.regular} cannot be applied for this class of solutions, since the condition that the vorticity be bounded is not enough to ensure the local-in-time  $H^s$ regularity of the velocity for $s>\frac{n}{2}+1$. The next result is established in \cite{Che2} for the two-dimensional case, where the effect of the singularity appears in the rate of convergence.

\begin{thm}\label{thm.bounded} Let $u^\nu_0=u_0$ be an $L^2$ perturbation of  a
smooth stationary solution to the Euler equations (see \cite[Definition 1.1]{Che2} for the precise definition).
If in addition $\omega_0={\rm curl}\, u_0 \in L^\infty (\mathbb{R}^2) \cap L^2 (\mathbb{R}^2)$, then
\begin{align*}
\| u^\nu  - u \|_{L^\infty (0,T; L^2)} \leq C \| \omega_0 \|_{L^\infty\cap L^2} (\nu T)^{\frac12 \exp (-C \|\omega_0 \|_{L^\infty \cap L^2} T)}\,.
\end{align*}
\end{thm}

Thus, although the inviscid limit is still verified in $L^\infty (0,T; L^2
(\mathbb{R}^2))$ even when the initial vorticity is merely in $L^\infty \cap
L^2$, the upper bound on the rate of convergence in Theorem \ref{thm.bounded}
decreases in time. This decrease, in fact, reflects the lack of Lipschitz
regularity for the velocity field of the Euler flow, e.g. $u\in L^\infty (0,T;
{\rm Lip} (\mathbb{R}^2))$.  On the other hand, the Lipschitz regularity of the
velocity is ensured for a class of vortex patches. Then,  as mentioned below,
the rate of convergence can be estimated uniformly in time for this class.
First, the definition of vortex patches is given as follows.

\begin{df}\label{def.patch} Let $n=2,3$ and $0<r<1$. A vector field $u\in L^p_\sigma (\mathbb{R}^n)$, $2<p<\infty$, is called a $C^r$ vortex patch if the vorticity $\omega ={\rm curl}\, u$ has the form
\begin{align}
\omega = \chi_A  \omega_i + \chi_{A^c} \omega_e\,,\label{def.patch.1}
\end{align}
where $A\subset \mathbb{R}^n$ is an open set of class $C^{1+r}$ and $\omega_i$, $\omega_e$ are compactly supported $C^r$ functions (vector fields when $n=3$). Here, $\chi_A$ and $\chi_{A^c}$ denote the characteristic functions of $A$ and $A^c=\mathbb{R}^n\setminus A$, respectively, and the condition $\omega_i \cdot {\bf n} = \omega_e \cdot {\bf n}$ is assumed on $\partial A$, which is always valid when $n=2$.
\end{df}

In Definition \ref{def.patch} the condition $\omega_i \cdot {\bf n} = \omega_e \cdot {\bf n}$ is assumed on $\partial A$ so that the divergence free condition ${\rm div}\, \omega=0$ is satisfied in the sense of distributions, which is necessary since $\omega={\rm curl}\, u$. The simplest vortex patch in the two-dimensional case is the constant vortex patch introduced in \cite{Maj}, where $\omega_i=1$, $\omega_e=0$, and $A$ is a bounded domain with $C^{1+r}$ boundary. In the three-dimensional case, the constant vortex patch cannot exist because of the requirement ${\rm div}\, \omega =0$. The classical reference for constant vortex patches is \cite{Che1}, where it is proved that if the initial vorticity is a $C^r$ constant vortex patch then it remains to be a $C^{r}$ constant vortex patch for all time; see also \cite{BeCo}.
Moreover, the velocity is bounded in $L^\infty (0,T; {\rm Lip}\, (\mathbb{R}^n))$ for all $T>0$.
The regularity of the vorticity field up to the boundary of general two-dimensional vortex patch is shown in \cite{De, Hua1}.
The result for the two dimensional case is extended to the three dimensional $C^{1+r}$ vortex patches by \cite{GaSa} but for a bounded time interval,
since, due to possible  vortex stretching mechanism, the existence of the Euler
solutions is only local in time.
The reader is also referred to \cite{Hua2} for a Lagrangian approach proving the regularity up to the boundary, and to \cite{Du} for the result in a bounded domain, rather than in the whole space.
Useful references about the study of vortex patches can be found in the recent paper \cite{Sue}.

The first result for the inviscid limit problem pertaining to vortex patches is
given in \cite{CoWu1},
where it is proved that the Navier-Stokes flows in $\mathbb{R}^2$ starting from a constant vortex patch converges to the constant vortex patch of the Euler flows with the same initial data.
In \cite{CoWu1} the convergence is shown in the topology of $L^\infty (0,T; L^2
(\mathbb{R}^2))$, $T>0$, for the velocity fields, with convergence rate $(\nu
T)^\frac12$.
The authors of  \cite{CoWu1} proved in \cite{CoWu2} the convergence of the vorticity fields in $L^p (\mathbb{R}^2)$, $2\leq p<\infty$.
The optimal rate of convergence is then achieved  for constant vortex patches in $\mathbb{R}^2$ by the results in \cite{AbDa}, as stated below.

\begin{thm}\label{thm.patch} Let $u^\nu_0=u_0$ be a constant $C^r$ vortex patch in $\mathbb{R}^2$.
Then there exists $C>0$ depending only on the initial data such that
\begin{align*}
\| u^\nu (t) - u(t) \|_{L^2}\leq C e^{C e^{Ct}} (\nu t)^\frac34 (1+\nu t)\,, \qquad t>0\,.
\end{align*}
\end{thm}
The rate $(\nu t)^\frac34$ is optimal in the sense that, if $u^\nu_0=u_0$ is a
circular vortex patch, i.e., $A$ is a disk, then an explicit computation leads
to the following bound from above and below (see \cite[Section 4]{AbDa}):
\begin{align*}
C (\nu t)^\frac34 \leq \| u^\nu (t) - u(t) \|_{L^2}\leq C' (\nu t)^\frac34\,, \qquad 0<\nu t\leq 1\,.
\end{align*}
The proof of Theorem \ref{thm.patch} in \cite{AbDa} is based on an energy method combined with the Littlewood-Paley decomposition, i.e., a dyadic decomposition in the Fourier variables.
In \cite{AbDa}, the authors used the bound  $u^\nu \in L^\infty (0,T; {\rm
Lip}\, (\mathbb{R}^2))$ uniformly in $\nu>0$, which was obtained in \cite{Dan1}
and extended in \cite{Dan2} for general $C^r$ vortex patches in $\mathbb{R}^n$,
to study the vanishing viscosity limit. The reader is also referred to
\cite{Hmi1} for the extension of \cite{Dan1}, and to \cite{Hmi2} for the
analysis when the constant patch has a singularity at its boundary.
Later it was pointed out in \cite{Masmoudi2007} that the uniform Lipschitz bound
of $u^\nu$ itself is not necessary and the optimal rate is proved only under a
Lipschitz bound of $u$.
More precisely,  in \cite[Theorems 3.2, 3.4]{Masmoudi2007}, the inviscid limit
is verified for general $C^r$ vortex patches in $\mathbb{R}^n$ in terms of Besov
spaces $B^s_{p,q}$, $s\in \mathbb{R}$, $1\leq p, q\leq \infty$, as follows.

\begin{thm}\label{thm.patch.general} Let $u_0^\nu=u_0$ be a $C^r$ vortex patch.
Assume that ${\rm curl}\, u_0\in \dot{B}^\alpha_{2,\infty} (\mathbb{R}^n)$,
$0<\alpha<1$, and that $u\in L^\infty (0,T; {\rm Lip}\, (\mathbb{R}^n))$. Then
any weak solution $u^\nu$  to the Navier-Stokes equations with initial data
$u^0$ satisfies the estimate
\begin{align*}
\| u^\nu (t) - u (t) \|_{L^2} \leq C (\nu t)^{\frac{1+\alpha}{2}}\,, \qquad 0<t<T\,.
\end{align*}
Here $C$ depends only on $u$ and $T$.
\end{thm}

\begin{rem}{\rm In Theorem \ref{thm.patch.general}, the class of weak solutions
is of Leray-Hopf type when $n=3$. When $n=2$ the solution is not always energy
finite due to the structure of the vortex patch, so the class of weak solutions
needs to be  suitably modified at this point. This is not an essential problem
here, for the global well-posedness of the Navier-Stokes equations in
$\mathbb{R}^2$ is known in various functional settings.
}
\end{rem}

As usual, the norms of the Besov spaces are defined in terms of the
Littlewood-Paley decomposition in the Fourier variables, cf.
\cite[Appendix]{Masmoudi2007}. The convergence rate in Theorem
\ref{thm.patch.general} depends on the exponent $\alpha$ for the regularity of
the vorticity, rather than the value $r$ for the regularity of the vortex patch.
An estimate of this type was  also obtained in \cite[Theorem 1.1]{AbDa} for the
two-dimensional case.
It is worthwhile to note that  the space $L^\infty (\mathbb{R}^n) \cap BV
(\mathbb{R}^n)$, where $BV$ denotes the space of functions of bounded
variation, is continuously embedded in $\dot{B}_{2,\infty}^{\frac12}
(\mathbb{R}^n)$ (see the proof of \cite[Lemma 4.2]{Masmoudi2007}).
Then,  a $C^r$ vortex patch $u$ satisfies ${\rm curl}\, u \in \dot{B}_{2,\infty}^{\frac12} (\mathbb{R}^n)$ if $\omega_i,\omega_e$ in Definition \ref{def.patch} belong to  $C^{\frac12} (\mathbb{R}^n)$ in addition.
For such a case, Theorem \ref{thm.patch.general} gives the convergence rate
$(\nu t)^\frac34$, which is known to be optimal in general.
The key observation here is that if $u\in L^\infty (0,T; {\rm Lip}\, (\mathbb{R}^n))$ and ${\rm curl}\, u_0 \in  \dot{B}_{2,\infty}^{\alpha} (\mathbb{R}^n)$, $\alpha\in (0,1)$, then the $\dot{B}_{2,\infty}^{\alpha} (\mathbb{R}^n)$ regularity is preserved under the convection by the flow $u$, and ${\rm curl}\, u\in L^\infty (0,T; \dot{B}_{2,\infty}^{\alpha} (\mathbb{R}^n))$ holds (for example, see the argument of \cite[Proposition 3.1]{AbDa}).
The formal proof of Theorem \ref{thm.patch.general} is given below, which
showcases the role of the regularity of ${\rm curl}\, u$ in estimating the
convergence rate, while a rigorous justification is given in
\cite{Masmoudi2007}. Since the difference $w=u^\nu-u$ solves the equations
\begin{align*}
\partial_t w - \nu \Delta w + u \cdot \nabla w + w\cdot \nabla u + w\cdot \nabla w + \nabla (p^\nu-p) = \nu \Delta u\,, \qquad {\rm div}\, w=0\,,
\end{align*}
a standard energy method, using  integration by parts, leads to
\begin{align*}
\frac12 \frac{\dd}{\dd t} \| w\|_{L^2}^2 & = - \nu \| \nabla w \|_{L^2}^2 - \langle w\cdot \nabla u, w\rangle _{L^2} - \nu \langle \nabla u, \nabla w \rangle _{L^2}\\
& \leq -\nu \| \nabla w\|_{L^2}^2 + \| \nabla u \|_{L^\infty} \| w\|_{L^2}^2 - \nu \langle \nabla u, \nabla w\rangle _{L^2}\,.
\end{align*}
Here $\langle \cdot, \cdot \rangle _{L^2}$ denotes the $L^2$ inner product.
The duality between $\dot{B}_{2,1}^{-\alpha}(\mathbb{R}^n)$ and $\dot{B}_{2,\infty}^\alpha (\mathbb{R}^n)$ implies
\begin{align*}
|\langle \nabla u, \nabla w\rangle _{L^2}| \leq C \| \nabla u \|_{\dot{B}^{\alpha}_{2,\infty}} \| \nabla w \|_{\dot{B}^{-\alpha}_{2,1}} = C \| \nabla u \|_{\dot{B}^{\alpha}_{2,\infty}} \|w \|_{\dot{B}^{1-\alpha}_{2,1}}\,,
\end{align*}
while ${\rm div}\, u=0$ and an  interpolation inequality (see \cite[Lemma 4.1]{Masmoudi2007}) yield:
\begin{align*}
\| \nabla u \|_{\dot{B}^{\alpha}_{2,\infty}}  = \| {\rm curl}\, u \|_{\dot{B}^{\alpha}_{2,\infty}}\,, \qquad \|w \|_{\dot{B}^{1-\alpha}_{2,1}}\leq C \| w\|_{L^2}^\alpha \| \nabla w\|_{L^2}^{1-\alpha}\,,
\end{align*}
respectively. Young's inequality finally gives:
\begin{align*}
\frac12 \frac{\dd}{\dd t} \| w\|_{L^2}^2  + \frac{\nu}{2} \| \nabla w\|_{L^2}^2 \leq \| \nabla u \|_{L^\infty} \| w\|_{L^2}^2 + C \nu \| {\rm curl}\, u \|_{\dot{B}_{2,\infty}^\alpha}^\frac{2}{1+\alpha} \| w\|_{L^2}^\frac{2\alpha}{1+\alpha}\,.
\end{align*}
Hence, in virtue of the bounds $u\in L^\infty (0,T; {\rm Lip}\, (\mathbb{R}^n))$
and ${\rm curl}\, u\in L^\infty (0,T; \dot{B}^\alpha_{2,\infty}
(\mathbb{R}^n))$, the estimate in Theorem \ref{thm.patch.general} follows from
 Gr\"onwall's inequality.

The results of Theorems \ref{thm.bounded} - \ref{thm.patch.general} show that the Euler flow describes the first order expansion of the Navier-Stokes flow  in viscosity in the energy norm for a suitable class of vortex patches.
It is then natural to investigate  higher-order expansions; however, this problem becomes highly nontrivial since the vorticity of a patch is discontinuous across the boundary of the patch.
Due to the smoothing effect of the viscosity term $\nu\Delta u^\nu$ in the Navier-Stokes equations, one has to introduce a fast scale in the higher order expansions which represents a viscous transition at the boundary of the vortex patch.
Hence, the boundary layer analysis comes to play an important role.
The first result in this direction is recently given in \cite{Sue}, where a
complete asymptotic expansion is provided in powers of $(\nu t)^\frac12$.

\vspace{0.3cm}

\noindent
{\bf (II) Singular solutions.} (II-2) Vortex sheets $\cdot$ Vortex filaments $\cdot$ Point vortices:
In many physical situations the vorticity of flows concentrates on a very small region.
Typical examples are vorticity fields called vortex sheets and vortex filaments or tubes,
and mathematically they are formulated as a class of flows the vorticities of which are Radon measures supported on hypersurfaces (vortex sheets) or curves (vortex filaments).

The velocity field associated to a vortex sheet satisfies the Euler equations on both side of the sheet and its tangential components jump across the sheet.  It is a natural model of flows with small viscosity after separation from rigid walls or corners. The simplest example in the two-dimensional case is the stationary flow
\begin{align}\label{ex.sheet}
u=(u_1,0)\,, \qquad
u_1 =
\begin{cases}
& -\frac12\,, \qquad x_2>0\,,\\
& \frac12\,, ~~\qquad x_2<0\,.
\end{cases}
\end{align}
The associated vorticity field is then given by $\omega (x) = \delta (x_2)$,
where $\delta (x_2)$ is the Dirac measure supported on $x_2=0$. However, this
simple exact solution is linearly unstable to small periodic disturbances, known
as the Kelvin-Helmholtz instability. More precisely, the linearization of the
Birkhoff-Rott equation, which is an equivalent formulation to the Euler
equations for vortex sheets under suitable conditions (\cite{LoNuScho}), around
\eqref{ex.sheet} has a solution growing exponentially like $e^{c|k|t}$ for each
$k$th Fourier mode of the interface corrugation (see e.g. \cite[Section
9.3]{MajdaBertozzi} for details). As a result, the solvability of the
Birkhoff-Rott equation is available only for analytic initial data
\cite{SuSuBarFri, CaO1,DuRo}, and it is known to be ill posed if one goes beyond
the analytic framework \cite{DuRo, CaO2}. As for the relation between the
analyticity and the regularity of the solution to the Birkhoff-Rott equations,
see also \cite{Le,Wu} and references therein.

The Navier-Stokes equations in $\mathbb{R}^n$ with vortex sheets as initial data are locally well-posed when $n=3$ and globally well-posed when $n=2$ since the velocity field in this class has enough local regularity to construct a unique solution. However, due to the underlying Kelvin-Helmholtz instability it is highly nontrivial to describe the behavior of the solutions in the inviscid limit. In \cite{CaSa06} the equations for viscous profiles of vortex sheets are presented in the two dimensional case when the radius of curvature of the sheet is much larger than the thickness of the layer, and those equations are solved within the analytic category.
The equations for viscous profiles in three dimensions are announced in \cite[Eq. (253)]{Sue}.
In both cases the equations exhibit the loss of one derivative.
A rigorous justification of the asymptotic expansion described in
\cite{CaSa06,Sue} is still lacking even for analytic data. At the same time, in
the two dimensional case, a weak solution to the Euler equations for vortex
sheet initial data can be constructed  with the help of  viscous solutions, if
the initial vorticity has a distinguished sign  (see \cite[Chapter
11]{MajdaBertozzi} for a detailed discussion). However, the verification of the
inviscid limit in generality as in the case of vortex patches seems to be out of
reach.

For a vortex filament, the vorticity field is concentrated on a curve in $\mathbb{R}^3$. The vortex filament is more singular than the vortex sheet.
The vorticity field of a vortex filament naturally belongs to  a Morrey space
$\mathcal{M}^{\frac32}(\mathbb{R}^3)$, which is an invariant space under the
scaling: $f_\lambda (x) = \lambda^2 f(\lambda x)$, $\lambda>0$ (see
\cite{GiMi}for the precise  definition of the space
$\mathcal{M}^{\frac32}(\mathbb{R}^3)$).
Vortex filaments possess an infinite energy in general.
Although the vortex filament can belong to the class of functional solutions to the Euler equations introduced by \cite{ChaDu}, due to the strong singularity and underlying vortex stretching mechanism in three dimensional flows, there seems to be no general existence result for vortex filaments as solutions to the Euler equations.
Alternatively, the self-induction equation (localized induction approximation) and its significant generalization have been used to understand the dynamics of vortex filaments; that dynamics is  out of the scope of this chapter and the interested reader is referred to \cite[Chapter 7]{MajdaBertozzi} and references therein.
On the other hand, the vortex filament belongs to an invariant space for the three-dimensional Navier-Stokes equations. Hence, a general theory is available for the unique existence of solutions to the Navier-Stokes equations with vortex filaments as initial data  under smallness condition on the scale-invariant norm \cite{GiMi}.

There are two typical vorticity distributions of vortex filaments: circular vortex rings with infinitesimal cross section, and exactly parallel and straight vortex filaments with no structural variation along the axis.
The first one, called vortex ring for simplicity here, corresponds to an axisymmetric flow without swirl.
Then the vorticity field of the vortex ring becomes a scalar quantity and is
expressed as a Dirac measure supported at a point  $(r,z)\in (0,\infty) \times
\mathbb{R}$ in cylindrical coordinates.
Recent results \cite{FeSve} show global existence of solutions to the
Navier-Stokes equations with vortex ring as initial data, and in fact their
uniqueness can be also proved in the class of axisymmetric flows without swirl;
a more detailed overview about this topic is available in another chapter of
this handbook, see \cite{GaMa}.
The inviscid limit problem for vortex rings is attempted in \cite{Mar3}, where
the initial vortex ring is slightly regularized depending on the viscosity and
the radius of the initial vortex ring is assumed to tend to infinity as the
viscosity goes to zero. The result of \cite{Mar3} is significantly extended in
\cite{BruMar}, where the radius of the initial vortex ring can be taken
independently of the viscosity.
A more precise description of the result in \cite{BruMar} is given as follows.
The velocity of the axisymmetric flow without swirl is expressed as
$u^{\nu}=(u_r^{\nu},0,u_z^{\nu})$ in cylindrical coordinates \ $(r,\phi,z)$, and
the vorticity field is identified with a scalar quantity,
$\omega^{\nu}=\partial_z u_r^{\nu} -\partial_r u_z^{\nu}$.
Then the evolution of $\omega^{\nu}$  in  cylindrical coordinates
is obtained from the equation
\begin{align}\label{eq.ring}
\partial_t \omega^{\nu} + (u_r^{\nu}\partial_r + u_z^{\nu}\partial_z ) \omega^{\nu} - \frac{u_r^{\nu}}{r} \omega^{\nu} = \nu \big ( \frac{1}{r}\partial_r (r \partial_r \omega^{\nu}) + \partial_z^2 \omega^{\nu} - \frac{\omega^{\nu}}{r^2} \big )\,.
\end{align}
Let
\begin{align*}
\Sigma_{(r_0,z_0)} (l) = \{ (r,z)\in (0,\infty)\times \mathbb{R}~|~ |r-r_0|^2 + |z-z_0|^2 <l^2 \}\,.
\end{align*}
Then a typical case of the result stated in \cite[Theorem 1.1]{BruMar} leads to the next theorem.
\begin{thm}\label{thm.filament} Assume that a sequence of initial vorticities $\displaystyle \{\omega_{0}^{\nu}\}_{0<\nu <\frac14}$ satisfies
\begin{align}\label{est.thm.filament.1}
{\rm supp}\, \omega_0^{\nu}  \subset \Sigma_{(1, 0)} (\nu^\frac12)\,, \qquad 0\leq \omega_0^{\nu} (r,z) \leq \frac{M}{\nu |\log \nu|}\,, \qquad  \int_{(0,\infty)\times \mathbb{R}} \omega_0^{\nu} \dd r \dd z = \frac{2a}{|\log \nu|}\,,
\end{align}
where $M>0$ and $a\in \mathbb{R}$ are constants independent of $\nu\in (0,\frac14)$.
Then there exists a sequence $\{(r_\nu (t), z_\nu (t))\}_{0<\nu<\frac14}$ in $(0,\infty) \times \mathbb{R}$ such that  for any $T>0$ and $f\in BC ([0,\infty)\times \mathbb{R})$,
\begin{align}\label{est.thm.filament.2}
\begin{split}
\lim_{\nu \rightarrow 0} r_\nu (t)  = 1\,, \qquad \lim_{\nu \rightarrow 0} z_\nu (t)  & = \frac{at}{4\pi}\,,\\
\lim_{\nu \rightarrow 0}  \bigg ( |\log \nu | \, \int_{\Sigma_{(r_\nu (t), z_\nu (t))}(D_\nu)} \omega^{\nu} (t) f  \dd r \dd z  \bigg ) & = 2 a f (1, \frac{at}{4\pi})\,,  \qquad 0<t\leq T\,, \\
 {\rm where}~~ D_\nu  & = C \nu^\frac12 \exp ( |\frac12 \log \nu|^\gamma)\,, \qquad 0<\gamma<1\,.
\end{split}
\end{align}
Here $\omega^{\nu}$ is the (unique) solution to \eqref{eq.ring} with initial
data $\omega_0^{\nu}$, and $C>0$ is a constant independent of $\nu$.
\end{thm}

Above BC stands for the space of bounded, continuous functions.
Theorem \ref{thm.filament} implies that, when the vorticity is initially sharply concentrated in an annulus, then
it remains concentrated during the motion even in the presence of small viscosity, and as $\nu$ goes to zero
the support of the vorticity evolves via a constant motion. The support
condition in \eqref{est.thm.filament.1} is related to the standard boundary
layer thickness, i.e., of order $\nu^\frac12$.
On the other hand, in \eqref{est.thm.filament.1} the vorticity is assumed to be logarithmically smaller than the standard scale.
This is related to the fact that the velocity field for the circular vortex ring
with infinitesimal cross section has a logarithmic singularity in the vertical
component around the location of the ring (that is, at $r=1$ in the setting of
Theorem \ref{thm.filament}). In virtue of the additional smallness of order
$\mathcal{O}(|\log \nu|^{-1})$ in \eqref{est.thm.filament.1} the vorticity is
translated with  finite speed in the inviscid limit, as described in
\eqref{est.thm.filament.2}.

Next, the other typical class of vortex filaments, i.e., parallel straight
vortex filaments, is discussed. By symmetry, these vortex filaments keep their
shape under the motion, and then the problem is reduced to the motion of the
points which are the intersection of the vortex filaments with the hyperplane
$\{x_3=0\}$. These are called point vortices, linear combinations of Dirac
measures in $\mathbb{R}^2$.
For the inviscid flow the motion of $N$ point vortices is formally described by the following Helmholtz-Kirchhoff system
\begin{align}\label{eq.HK.system}
\frac{\dd }{\dd t} z_i (t) = \frac{1}{2\pi} \sum_{j\ne i} \alpha_j \frac{(z_i (t) - z_j (t))^\bot}{|z_i (t) - z_j (t)|^2}\,, \qquad z_i (0) = x_i\,.
\end{align}
Here $i,j=1,\cdots,N$ and each $x_i$ denotes the initial location of the $i$th point vortex with circulation $\alpha_i\in \mathbb{R}$.
It is possible to realize the point vortices as  functional solutions, introduced in \cite{ChaDu}, to the two-dimensional Euler equations. A further relation of point vortices with the Euler equations is shown by \cite{MarPu}. On the other hand, the two-dimensional Navier-Stokes equations are known to be globally well-posed when the initial vorticity field is given by the point vortices of the form \eqref{est.thm.pv.1.1} below. In particular, the uniqueness of solutions is also available, which is proved in \cite{GiMiO,Ka3} under smallness condition on the total variation $\sum_{i=}^N |\alpha_i|$, and in  \cite{GaGa} without any smallness condition on the size of $\alpha_i$.

The inviscid limit problem for point vortices is rigorously analyzed in  \cite{Mar1,Mar2, Gal} in the time interval in which the Helmholtz-Kirchhoff system is well-posed and vortex collisions do not occur. The next result is due to \cite[Theorem 2]{Gal}.
\begin{thm}\label{thm.pv.1} Assume that the point vortex system \eqref{eq.HK.system} is well-posed on the time interval $[0,T]$. If the initial vorticity field is given by
\begin{align}\label{est.thm.pv.1.1}
{\rm curl}\, u_0 = \sum_{i=1}^N \alpha_i \, \delta (\cdot - x_i)\,,
\end{align}
then the vorticity field $\omega^{\nu}={\rm curl}\, u^{\nu}$ of the solution to the Navier-Stokes equations \eqref{e.NS} converges to $\displaystyle \sum_{i=1}^N \alpha_i \delta (\cdot - z_i (t))$ as $\nu \rightarrow 0$ in the sense of measures for all $t\in [0,T]$,
where $z(t) = (z_1(t), \cdots, z_N (t))$ is the solution of \eqref{eq.HK.system}.
\end{thm}
Theorem \ref{thm.pv.1} shows that the distribution of the vorticity field of the Navier-Stokes flow starting from \eqref{est.thm.pv.1.1} is described by the point vortex system \eqref{eq.HK.system} in the inviscid limit.
A result similar to  Theorem \ref{thm.pv.1} was first obtained in
\cite{Mar1,Mar2}, where the initial point vortices are slightly regularized
depending on viscosity.
Theorem \ref{thm.pv.1} provides  information on the location of viscous
vortices,
while little information is  available about the shape of the viscous vortices, i.e., about the viscous profile of each point vortex. Since the vorticity field $\omega^{\nu}$ obeys the nonlinear heat-convection equations
\begin{align*}
\partial_t\omega^{\nu} + u^{\nu} \cdot \nabla \omega^{\nu} =\nu \Delta \omega^{\nu} \,, \qquad \omega^{\nu}|_{t=0} = \sum_{i=1}^N \alpha_i \delta (\cdot - x_i)\,,
\end{align*}
once $\omega^{\nu}$ is constructed by solving the above system, it is then decomposed as
\begin{align}\label{eq.decompose}
\omega^{\nu} = \sum_{i=1}^N \omega_i^{\nu}\,,
\end{align}
where $\omega_i^{\nu}$ is the solution to the heat-convection equations with
initial data which is a Dirac measure supported at $x_i$:
\begin{align}\label{eq.omega_i}
\partial_t\omega_i^{\nu} + u^{\nu} \cdot \nabla \omega_i^{\nu} =\nu \Delta \omega_i^{\nu} \,, \qquad \omega_i^{\nu}|_{t=0} = \alpha_i \delta (\cdot - x_i)\,.
\end{align}
An interesting and important question is then to determine the correct
asymptotic profile of each $\omega_i^{\nu}$.
It is worthwhile to recall here that, if the initial vorticity field is a Dirac
measure supported at the origin, $\alpha \delta(0)$, then the unique solution to
\eqref{e.NS} is explicitly given by the Lamb-Oseen vortex with circulation
$\alpha$  (\cite{GaWa}):
\begin{align}\label{def.Lamb-Oseen}
\alpha U^{\nu} (t,x) = \alpha \frac{x^\bot}{2\pi |x|^2} ( 1- e^{-\frac{|x|^2}{4\nu t}})\,, \quad \qquad \big ({\rm curl}\, \alpha U^{\nu}  \big ) (t,x) = \frac{\alpha}{\nu t} G (\frac{x}{\sqrt{\nu t}})\,, \qquad G (x) = \frac{1}{4\pi} e^{-\frac{|x|^2}{4}}\,.
\end{align}
Hence, the function $\alpha_i G$, a Gaussian with  total mass $\alpha_i$,  is a
natural candidate for the viscous profile of each $\omega_i^{\nu}$.
The next result (\cite[Theorem 3]{Gal}) establishes the asymptotic
expansion of $\omega_i^{\nu}$ in the limit $\nu\rightarrow 0$, and shows that
the asymptotic profile of  each $\omega_i^{\nu}$ is indeed given by a
two-dimensional Gaussian.

\begin{thm}\label{thm.pv.2} Assume that the point vortex system \eqref{eq.HK.system} is well-posed on the time
interval $[0,T]$, and let $\omega^{\nu}={\rm curl}\, u^{\nu}$ be the vorticity
field of the solution to \eqref{e.NS} with initial data the vorticity of which
is given by \eqref{est.thm.pv.1.1}.
If $\omega^{\nu}$ is decomposed as in \eqref{eq.decompose}, then the rescaled
profiles $w_i^{\nu}$, defined by
\begin{align}\label{def.scale}
\omega_i^{\nu} (t,x) = \frac{\alpha_i}{\nu t} w_i^{\nu} (t, \frac{x-z_i^\nu
(t)}{\sqrt{\nu t}}),
\end{align}
satisfy the estimate
\begin{align}\label{est.thm.pv.2.1}
\max_{i=1,\cdots, N} \| w_i^{\nu} (t) - G\|_{L^1}\leq C\frac{\nu t}{d^2}, \qquad
 \qquad t\in (0,T]\,,
\end{align}
where $d=\displaystyle \min_{t\in [0,T]}\min_{i\ne j}| z_i (t) - z_j (t)|>0$. Here $z_i^\nu$ is the solution to
a regularized point vortex system (see \cite[Eq. (19)]{Gal}) .
\end{thm}

In \cite{Gal}, the estimate \eqref{est.thm.pv.2.1} is proved in a stronger
topology, that of  a weighted $L^2$-space that is embedded in $L^1
(\mathbb{R}^2)$. The behavior of $z_i^{\nu}$ in \eqref{def.scale} is well
approximated by $z_i$ in the limit $\nu\rightarrow 0$ with an exponential order
\cite[Lemma 2]{Gal}, and hence, Theorem \ref{thm.pv.2} verifies the asymptotic
expansion of the viscous vortices around $z_i(t)$, $t\in (0,T]$.
In order to show Theorems \ref{thm.pv.1} and \ref{thm.pv.2}, one has to take into account the interaction between the viscous vortices in the inviscid limit. In particular, for each viscous vortex around $z_i$ the velocity fields produced by the other vortices play a role in  the background flow. The interactions results in  a deformation of each viscous vortex. The analysis of this interaction is the key to prove Theorems \ref{thm.pv.1} and \ref{thm.pv.2}, and in fact it requires a detailed investigation of  higher-order expansions due to the strong singularity of the flows. This approach is validated in a rigorous fashion in \cite[Theorem 4]{Gal}. One consequence is  that each viscous vortex is deformed elliptically through the interaction with the other vortices.

\section{Inviscid limit problem with physical boundary}\label{sec.phys.b}

This section is devoted to the analysis of a viscous fluid at very low
viscosity moving in a domain $\Omega$ with physical boundaries.
The behavior of the fluid is markedly influenced by the types of boundary conditions imposed.
The case when rigid walls may move only parallel to itself,
as in Taylor-Couette flows, will be the focus in this section. In this case, the fluid domain (assumed to be
smooth) is fixed and both the viscous, inviscid flows must satisfy the {\em
no-penetration} condition at the boundary:
\begin{equation} \label{e.nopenBC}
     u\cdot {\bf n}=0,
\end{equation}
where $u$ is the fluid velocity and ${\bf n}$ is the unit outer normal to the domain, respectively.
For ideal fluids, the no-penetration condition is the only one that can be
imposed on the flow. It is often referred to, somewhat incorrectly, as a {\em
slip boundary condition}, because the fluid is allowed to slip, but no slip
parameters are specified. For viscous fluids, there are several possible
boundary conditions that are physically consistent. The simplest, and most
difficult one from the point of view of the vanishing viscosity limit, is the
{\em no-slip} boundary condition, where \eqref{e.nopenBC} is
complemented with
the following condition on the tangential fluid velocity at the boundary:
\begin{equation} \label{e.noslipBC}
   u_{\text{tan}} = V,
\end{equation}
where  $V$ is the velocity of the boundary and $u_{\text{tan}}$ is the tangential
component of the velocity. This is the boundary condition originally proposed by
Stokes.
If the friction force is prescribed at the boundary,
then one obtains {\em Navier friction} boundary condition, which allow for slip
to occur:
\begin{equation} \label{e.slipBC}
  u\cdot {\bf n}=0, \qquad  [S(u)\, {\bf n} + \alpha u]_{\text{tan}} = 0,
\end{equation}
where $\alpha\geq 0$ is the friction coefficient and $S(u)$ is the viscous
stress
tensor, which coincides for Newtonian fluids with a multiple of the rate of
strain tensor $(\nabla u +\nabla u^T)/2$. Above,
given a vector field $v$ on the boundary of $\Omega$,
$v_{\text{tan}}$ means the component of $v$ tangent to
$\partial\Omega$, and if $M$ is a
matrix, $M\,v$ denotes matrix-vector multiplication. These
boundary conditions are the ones originally proposed by
Navier and derived by Maxwell in the context of gas dynamics.
In absence of friction ($\alpha=0$), the Navier boundary
condition reduces to the condition that the tangential component of the shear stress  be zero at the boundary.
They are also called {\em stress-free} boundary conditions in the literature.

The term $\alpha u$ can be replaced by $A u$, where $A$ is a
(symmetric) operator. In particular, if $A$ is the shape operator on the
boundary of the domain, then the generalized Navier boundary conditions reduce
to the following \cite{BdVC12,GK12}:
\begin{equation} \label{e.freeBC}
  u\cdot {\bf n}=0, \qquad \curl \, u \times {\bf n} =0,
\end{equation}
sometimes called {\em slip-without-friction} boundary conditions, as they
reduce to \eqref{e.slipBC} with $\alpha=0$ on flat portions on the boundary. In
two space dimensions, they are also referred to as {\em free} boundary
conditions, given that the second condition reduces to \ $\curl \, u = 0$
\cite{Lions69}, that is, there is no vorticity production at the boundary. In
higher dimensions, free boundary conditions lead to an overdetermined system,
but are still compatible with the time evolution of the fluid since the boundary
is characteristic if the initial data satisfies the same condition.
Other types of slip boundary conditions based on vorticity can be imposed
\cite{BNP04}, but they have been studied less in the literature.

Impermeability of the boundary implies that boundary conditions need to be
imposed on a characteristic boundary, which complicates the analysis further.
If the boundary is not characteristic, as is the case of a permeable boundary
with injection and suction, where the normal velocity at the boundary is
prescribed and non-zero, the zero-viscosity limit holds at least for short
time \cite{TW02}, as boundary layers can be shown to be very weak (of
exponential type). The case of a non-characteristic boundary is discussed in
Section \ref{s.nonchar}.

As already mentioned, one of the main obstructions to establishing the vanishing viscosity limit in
the presence of boundaries is the formation of a viscous boundary layer, where
the behavior of the flow cannot be approximated by that of an inviscid flow. A
formal asymptotic analysis using $\sqrt{\nu}$ as small parameter leads to a
reduced set of effective equations for the leading order velocity term in the
boundary layer, the so-called Prandtl equations (we refer to \cite{Schlichting}
for a historical perspective). Significant advances have been made in the analysis
of these equations, which exhibit instabilities, possible blow up, and
ill posedness. The analysis on the Prandtl equations and the connection between
solutions to Prandtl equations and the validity of the vanishing viscosity
limit will be reviewed in Section \ref{s.Prandtl}.

The discussion here will be confined to the classical case of the vanishing viscosity
limit for incompressible, Newtonian fluids, although some partial results are
available in the important cases of compressible flows
\cite{XY99,WX05,S14}, MHD \cite{XXW09,XL11}, convection in porous media
\cite{KTW11}, and for non-Newtonian (second-grade) fluids \cite{BILN12,LMNTZ15}.
In the ensuing discussion, the case of unsteady flows
in bounded, simply connected domains or a half space will be the main focus. The interesting case of
exterior or multiply-connected domains, such as flow outside one or more
obstacles, brings in additional difficulties, for example the infinite energy in
the vorticity-velocity formulation of the fluid equations in two dimensions (see
\cite{ILN09,KLN09,S12} and references therein).

The introduction of slip makes the vanishing viscosity limit more tractable,
essentially because the viscous boundary layer is weak compared to the outer
Euler solution and it is possible to obtain {\em a priori\/} bounds on higher norms
that are uniform in viscosity. The discussion starts, therefore, with
slip-type boundary conditions in Section \ref{s.slip} and continues with the
more challenging case of no-slip boundary conditions in Section \ref{s.noslip}.

In the remainder of this section,
the solution to the Navier-Stokes equations \eqref{e.NS} is denoted by $u^\nu$,
and the solution to the Euler equations \eqref{e.EE} is written as $u^0$.
Unless otherwise stated, it is also assumed that the initial data for the
Navier-Stokes equations is {\em ill prepared}, that is, the initial velocity
$u^\nu(0)$ is divergence free and
satisfies the no-penetration condition at the boundary, but not necessarily
viscous-type, e.g. friction or no-slip, boundary conditions.
This assumption will allow, in particular, to take the same initial data for \eqref{e.NS}
and \eqref{e.EE}:
\[
         u^\nu(0)=u^0(0)=u_0,
\]
although this assumption can, and will at times, be relaxed.
When the data is ill prepared, there is a corner-type singularity at $t=0$, $x\in \pa\Omega$ with two types of layers for the viscous evolution, an {\em initial layer} and the boundary layer. The initial layer is particularly relevant when the limit Euler flow is steady, as it affects vorticity production in the limit.

\subsection{Case of slip-type boundary condition}\label{s.slip}

If the viscous boundary layer is of size much smaller than predicted by Prandtl
asymptotic theory, one expects that the fluid will effectively slip at the
boundary (see \cite{LBS07} for a review of experimental results).  This
situation is more likely the rougher the boundary is. In fact,
it can be shown rigorously in certain situations that
homogenization of the no-slip boundary condition on a highly oscillating
boundary will give rise in the limit to slip boundary conditions
(see \cite{JM01,BFN10,MNN13} and references therein).

With Navier friction boundary conditions, the vanishing viscosity limit holds in
two and three space dimensions under additional regularity conditions on the
initial data, even if vorticity is produced at the boundary and the boundary is
characteristic. In two space dimensions, the problem can be studied in the vorticity-velocity
formulation, as the Navier friction condition gives rise to a useful boundary
condition for vorticity, namely:
\begin{equation} \label{e.vorticityslipbc}
  \omega_\nu= (2 \kappa - \alpha) u^\nu_{\text{tan}},
\end{equation}
where $\kappa$ is the curvature of the boundary. The $2$D Navier-Stokes
initial value problem in vorticity-velocity formulation then reads:
\begin{equation} \label{e.NSvort2D}
 \begin{cases}
  \omega^\nu_t + u^\nu \cdot \nabla \omega^\nu = \nu \Delta \omega^\nu, &
\text{ on } (0,T)\times \Omega,\\
   u^\nu = K_\Omega[\omega^\nu], & \text{ on } (0,T)\times \Omega,\\
   \omega^\nu(0)=\omega_0, & \text{ on }  \Omega \\
   \omega_\nu= (2 \kappa - \alpha) u^\nu_{\text{tan}},  \text{ on } (0,T)\times
   \partial\Omega,
\end{cases}
\end{equation}
with $K_\Omega$ the Biot-Savart operator associated to the domain $\Omega$.
A bound on the initial vorticity in $L^p$, $p>1$,
ensures, in particular, a uniform-in-time bound on the $L^p$-norm of
the vorticity and hence global existence and uniqueness of a strong solution to
\eqref{e.NSvort2D}.
In general, a solution of this problem is not a weak solution of
\eqref{e.NS}, even if $\omega_0\in L^p(\Omega)\cap L^1(\Omega)$ for
some $p>1$, as $u^\nu$ is not of finite energy. However, this is the case if
$\Omega$ is a bounded domain.

The zero-viscosity limit was first established for bounded initial vorticities
and forcing in \cite{CMR98} and for unbounded forcing in
\cite{Ru06}. This result was then extended to initial vorticities in $L^p$,
$p>2$, \cite{LNP05}, and to initial vorticities in the Yudovich uniqueness
class, \cite{K06}, when the forcing is zero.
To be precise, the result in \cite{LNP05} is stated below, which applies to the
most general class of initial data.
\begin{thm} \label{thm.2DNavierConv}
Let $\omega_0\in L^p(\Omega)$, $p>2$, and let $u_0=K_\Omega[\omega_0]$. Let
$u^\nu$ be the unique weak solution of \eqref{e.NS} with $u^\nu(0)=u_0$.
Then there exists a sequence $\nu_k\to 0$ and a distributional solution $u^0$ of
the Euler equations \eqref{e.EE} with initial data $u_0$, such that $u^\nu_k\to
u^0$ strongly in
$C([0,T];L^2(\Omega))$ as $k\to \infty$.
\end{thm}
The proof relies on {\em a priori} bounds on higher Sobolev norms for the
velocity that are uniform in viscosity. These bounds in turn allow by
compactness, via an Aubin-Lions-type lemma,  to pass to a limit along
subsequences as $\nu\to 0$. Because compactness arguments are used, the
proof does not give rates of convergence of $u^\nu$ to $u^0$. Uniqueness of the
limit can be guaranteed {\em a posteriori} if $\omega_0$ is sufficiently
regular. The key ingredient in establishing
the vanishing viscosity limit is an a priori bound on the vorticity in
$L^\infty([0.T], L^p)$ uniform in viscosity (for e.g. $\nu \in (0,1]$). This
estimate in turn is achieved through the use of the maximum principle, which
requires estimating the $L^\infty$-norm  of the velocity at the boundary,
except in the case of free boundary conditions.
The condition $p>2$ allow to obtain such an estimate via the Sobolev embedding
theorem, but it is not expected to be sharp. A more natural condition would be
$p>1$, which would ensure passing to a limit in non-linear terms.
 When $\Omega=\R^2_+$ the zero-viscosity limit in the energy norm as
in Theorem \ref{thm.2DNavierConv} is justified
by \cite[Theorem 2]{Pa14} even in the case when the friction coefficient $\alpha$ depends on the viscosity
as $\alpha = \alpha' \nu^{-\beta}$ with some constants $\beta\in [0,1)$ and $\alpha' > 0$.
Interestingly, the instability of the Prandtl boundary layer shown by \cite{Grenier00} in the no-slip case (see Theorem \ref{thm.Grenier} below) can be proved as well in this {\it viscosity-dependent} slip condition when $\beta =\frac12$ (\cite[Theorem 3]{Pa14}),
for which the $L^2$ convergence for velocity is valid.
This result implies that there is an essential discrepancy between the validity
of the zero-viscosity limit in the energy norm and the validity of the Prandtl
boundary layer expansion.

In three space dimensions, it is generally only possible to establish the
vanishing viscosity limit for more regular Euler initial data, namely $u_0$
in the Sobolev space $H^s$, $s>\frac52$, and only for the time of existence of
the strong Euler solution. Weak {\em wild} solutions are generically too
irregular to allow passage to the vanishing viscosity limit even in the full
space (this point will be revisited in Section \ref{s.noslip}).
For regular initial data $u_0$, the results in Theorem \ref{thm.2DNavierConv}
extend to three space dimensions \cite{IP06}. In fact, since the maximum
principle for the vorticity no longer holds due to vortex stretching, a direct
energy estimate on the velocity is employed and, consequently, the
approximating sequence of Navier-Stokes solutions $u^{\nu_k}$ can be taken
in the Leray-Hopf class of weak solutions.

Uniform bounds in viscosity in higher Sobolev norms $H^s$,
$s>\frac52$, do not hold, due to the presence of a boundary layer, except in the
case of free boundary conditions. Nevertheless, it is possible to
show convergence of the Navier-Stokes velocity $u^\nu$ to the Euler solution
$u^0$ uniformly in space and time, but utilizing so-called co-normal Sobolev
spaces \cite{MR12}. These are Sobolev spaces that, at the boundary,
give control on tangential derivatives, defined as
\begin{equation}\label{e.conormalSpace}
  H^m_{\text{co}} (\Omega) =\{ f\in L^2(\Omega), Z^\alpha\,f\in L^2(\Omega),
|\alpha|\leq m \},
\end{equation}
where $Z^\alpha$ is a vector-valued differential operator which is tangent to
$\partial\Omega$. Such co-normal spaces will be employed in studying the no-slip
case as well in Section \ref{s.noslip}.
Because of the Navier boundary condition, it is possible to obtain a
bound in $L^\infty((0,T)\times \Omega)$ uniform in viscosity on the full
gradient of the Navier-Stokes velocity using co-normal spaces of high enough
regularity. Then, this Lipschitz control on the Navier-Stokes solution allows to
pass to the limit as viscosity vanishes by using compactness again.
Set
\[
      V^m=\{ f \in H^m_{\text{co}} (\Omega), \nabla f  \in H^m_{\text{co}} (\Omega),
         {\rm div}\, f=0 \}
\]
Then the convergence result in \cite{MR12} is stated as follows.
\begin{thm}\label{thm.3DNavierConv}
Fix $m\in \mathbb{Z}_+$, $m>6$. Let $u_0\in V^m$.
Assume in addition that $\nabla u_0\in W^{1,\infty}_\text{co}$. Let $u^\nu\in
C([0,T],V^m)$,  be the strong solution of the Navier-Stokes equations
\eqref{e.NS} with initial data $u_0$ and boundary conditions \eqref{e.slipBC}.
Then, there exists a unique solution to the Euler equations \eqref{e.EE} with
initial data $u_0$ and boundary condition \eqref{e.nopenBC}, $u^0\in
L^\infty((0,T), V^m)$, such that $\nabla u^0\in L^\infty((0,T),W^{1,\infty})$
and such that
\[
    \|u^\nu - u^0\|_{L^\infty((0,T)\times\Omega)} \to 0, \qquad \text{as } \nu
\to 0.
\]
\end{thm}
Note that in Theorem \ref{thm.3DNavierConv} the high regularity of the initial data is needed to
control the pressure at the boundary.

One can interpret the uniform convergence of Theorem \ref{thm.3DNavierConv} in
terms of boundary layer analysis.
For simplicity the boundary is assumed to be flat, or the domain $\Omega$ is locally identified with the half space $\mathbb{R}^n_+$ by writing a point $x\in \Omega$ as $x=(x',z)$ with $x'\in \mathbb{R}^{n-1}$,
where $n=2,3$ is the space dimension again, and $z>0$.
Thus, the boundary of $\Omega$ is identified locally with $z=0$.
Then, the following asymptotic expansion of the viscous velocity holds in the energy space
$L^\infty([0,T],L^2(\Omega))$ \cite{IS11}:
\begin{equation} \label{e.NavierPrandtlExp}
   u^\nu (t,x) = u^0(t,x) + \sqrt{\nu} U(t, x', \frac{z}{\sqrt{\nu}}) + O(\nu),
\end{equation}
where $U$ is a smooth, rapidly decreasing boundary layer profile on
$\mathbb{R}^n_+$.
Hence, the boundary layer has the same width as predicted by the Prandtl theory for the case
of no-slip boundary condition, but small amplitude. By contrast, the amplitude
of the boundary layer corrector to the Euler velocity can be of order one if
no-slip boundary conditions are imposed.
The validity of the expansion \eqref{e.NavierPrandtlExp} implies, in
particular, that in general one cannot expect strong convergence in high
Sobolev norm, since then, by the trace theorem, the limit Euler solution would
satisfy the Navier-slip boundary condition (the so-called {\em strong
zero-viscosity limit}). For the case of slip without friction, it was shown in
\cite{BdVC12} that the boundary condition \eqref{e.freeBC} is not necessarily
preserved under the Euler evolution in three space dimensions if the boundary is
not flat. The strong zero-viscosity limit does hold for free boundary
conditions, hence for arbitrary smooth domains in two   space
dimensions, and for slip-without friction boundary conditions on domains with
flat boundary \cite{BdVC10,BdVC13,XiaoXin07,BdVC14}. In fact, only the initial
data need to satisfy the stronger free condition at the boundary \cite{BS12}.
It is interesting to note that the main difficulty in dealing with non-flat
boundaries comes from the non-vanishing of certain integrals in the energy
estimate due to the convective term $u^\nu\cdot \nabla u^\nu$. In the no-slip
case, the interaction between convection and the boundary layer is thought to
be a main obstruction to the validity of the zero-viscosity limit.

The last part of this subsection is about an approach to the zero viscosity limit
that yields rate of convergence in viscosity. The point of departure is the
expansion \eqref{e.NavierPrandtlExp}. If this expansion is valid, then $u^0 +
\sqrt{\nu}V$ approximates the Navier-Stokes solution. One can,
therefore,
define an approximate Navier-Stokes solution $u^\nu_{\text{approx}}$ in terms of
an {\em
outer solution\/} $u_{\text{ou}}^\nu$ valid away from the walls, and an {\em
inner solution\/}
$u_{\text{in}}^\nu$ valid near the walls (cf. \cite{VL57,Lions73}). The parabolic
nature of the Navier-Stokes
equations suggests that the viscous effects are felt in a thin layer close to
the boundary of width $\sqrt{\nu}$ (this idea goes back to the original work of
Prandtl in fact, again the reader is referred to \cite{Schlichting} for more details).
For simplicity $\Omega$ is assumed to be a half-space, as in \eqref{e.NavierPrandtlExp}.
To define the inner solution on a fixed domain independent of viscosity, it is
convenient to introduce the stretched variable $Z=z/\sqrt{\nu}$. If the
zero-viscosity limit holds, one then expects that a regular asymptotic expansion
for $u^\nu_{\text{in}}$ is valid, that is:
\[
      u^\nu_{\text{in}}(t,x)= \sum_{k=0}^\infty \nu^{\frac{k+1}{2}}\,
\theta^\nu_k(t,x', Z),
\]
where $\theta^\nu_k$ is the $k$th-order {\em corrector} to the outer flow.
It should be stressed that the correctors are not assumed to be independent of
viscosity and their amplitude is dictated by the equations of motions and by
the boundary conditions.
Similarly, the outer solution should have a regular expansion of the form:
\[
      u^\nu_{\text{ou}}(t,x)= \sum_{k=0}^\infty \nu^{\frac{k}{2}}\,
 u^\nu_k(t,x',z),
\]
with $u^\nu_0$ independent of $\nu$. In fact, $u^\nu_0=u^0$, the Euler
solution. Consistency of the formal asymptotic expansion gives effective
equations for the flow correctors, together with boundary and initial
conditions, from the Navier-Stokes and Euler equations.
The goal is then to derive the regularity of the correctors and their decay away from the boundary $Z=0$ from
the effective equations. The regularity and decay properties typically depend
on compatibility conditions between the initial and boundary data. Using these
properties, norm bounds on the error $u^\nu-u^\nu_{\text{approx}}$ can then be
obtained from the Navier-Stokes equations via energy estimates.

This approach was used in \cite{GK12} to establish the zero-viscosity limit and
rates of convergence under generalized Navier boundary conditions on any
smooth, bounded domain; see Theorem \ref{thm.BLNavier} below.
Both \cite{IS11} and \cite{GK12} employ {\em linear} boundary correctors,
another measure of the weakness of the boundary layer
under Navier conditions, but the corrector in \cite{GK12} is constructed using
a covariant formulation and is coordinate independent. Under geodesic
boundary normal coordinates in a tubular neighborhood of the boundary, it has an
explicit form, which allows to prove uniform space-time bounds on the error
$u^\nu-u^0$ even close to the boundary.
\begin{thm}\label{thm.BLNavier}
Denote by $\Gamma_a$ a tubular neighborhood of $\partial\Omega$ interior to
$\Omega$ of width $a>0$. Fix $m>6$ and assume $u_0\in H^m(\Omega)$.
Let $u^\nu$ be the unique, strong solution of \eqref{e.NS}  with
generalized Navier boundary conditions and initial data $u_0$.
Let $u^0$ be the unique strong solution of \eqref{e.EE} with the no-penetration
boundary conditions and initial data $u_0$.  Then:
\[
  \begin{aligned}
    &\|u^\nu-u^0\|_{L^\infty([0,T],L^2(\Omega))} \leq \kappa
     \, \nu^{\frac{3}{4}},
     &\|u^\nu-u^0\|_{L^\infty([0,T],H^1(\Omega))} \leq \kappa
     \, \nu^{\frac{1}{4}} \\
     &\|u^\nu-u^0\|_{L^\infty([0,T] \times \Gamma_a))} \leq \kappa \,
    \nu^{\frac{3}{8} - \frac{3}{8(m-1)}},
      &\|u^\nu-u^0\|_{L^\infty([0,T] \times \Omega\setminus \Gamma_a))} \leq
    \kappa \, \nu^{\frac{3}{4} - \frac{9}{8m}},
  \end{aligned}
\]
where $\kappa$ is a constant independent of $\nu$.
\end{thm}
In two space dimensions, under free boundary conditions, the use of correctors
allows to study the boundary layer for the vorticity and obtain rates of
convergence in Sobolev spaces \cite{GJ13}.
Lastly, in the context of Navier-type boundary conditions, the zero-viscosity limit has been used as a mean to establish existence of solutions to the Euler equations (see e.g. (see  (\cite{Y1963},
\cite[pp.~87--98]{Lions69}, \cite{Bardos1972}, and  \cite[pp.~129--131]{L1996}). Whether it is also a selection mechanism for uniqueness of weak solutions in two space dimensions remains open.

\subsection{Case of no-slip boundary condition}\label{s.noslip}

The classical case of no-slip boundary conditions \eqref{e.noslipBC} is perhaps
the most relevant in applications, and the most challenging to study from a
mathematical point of view. The main difficulty stems from the formation of a
possibly strong boundary  layer (of amplitude order one in viscosity) in flows
at sufficiently high Reynolds numbers. It is experimentally observed (see e.g.
the classical experiments of flow around a solid sphere \cite{VanDyke82}) that
laminar boundary layers, where the flow lines are approximately parallel to the
boundary, destabilize and detach from the boundary, a phenomenon known as
{\em boundary layer separation}.
This layer separation is due to the presence of an adverse pressure gradient
that leads to stagnation first and then flow reversal in the layer. In the
unsteady case, the connection between the vanishing viscosity limit and the
stability of the boundary layer has still not been completely clarified. In
particular, there are no known analytical examples of  unsteady flows
where layer separation occurs. Nevertheless, a connection can be made in terms
of vorticity production at the boundary. The mismatch between no penetration and
no slip at the boundary leads potentially to the creation of large gradients of
velocity in the layer, in particular normal derivatives of tangential components
of the velocity at the boundary.
While the creation of a boundary layer can be a purely diffusive effect, it is its interaction with strong inertial terms that is thought to lead to boundary layer separation. Therefore, one expects that in the context of the Navier-Stokes equations linearized around a non-trivial profile, i.e., Oseen-type equations, it should be possible to establish the zero-viscosity limit. This is indeed the case at least if the Oseen profile is regular enough and under some compatibility conditions between the initial and boundary data \cite{Aleks83,Aleks86,TWOseen96,TWOseen97,LSOseen01} (see also \cite{Gie13} for incompatible data for the Stokes equation).

Whether the vanishing viscosity limit holds generically even for short time under no-slip boundary conditions is largely an open problem.
An asymptotic ansatz for the velocity similar to  \eqref{e.NavierPrandtlExp}, but taking into account the amplitude of the boundary layer profile, that is, postulating an expansion for the velocity of the form:
\begin{equation} \label{e.noslipPrandtlExp}
 u^\nu(t,x) = u^0(t,x) + \theta(t,x',\frac{z}{\sqrt{\nu}}) + O(\nu^{\frac12})\,,
\end{equation}
leads to the classical Prandtl equation, which will be discussed in more details in Section \ref{s.Prandtl}.
One brief remark here is that the Prandtl equations are well posed only under strong conditions on the flow, such as when boundary and the data have some degree of analyticity \cite{Asano88,SC98,LCS03,CLM13,KV13} or the data is monotonic in the normal direction to the boundary \cite{Oleinik66,OleinikSamokhin99,KMVW14}.
The most classical result verifying \eqref{e.noslipPrandtlExp} is \cite{SC98} in the analytic functional framework, after the pioneering work of \cite{Asano88,Asano88'}.
The result of \cite{SC98} is stated here only in an intuitive manner without introducing the precise definition of function spaces.
\begin{thm}\label{thm.SC98} Suppose that given data and the Euler flow are analytic in all variables $x=(x',z)$.
Then the Prandtl asymptotic expansion \eqref{e.noslipPrandtlExp} is valid for a short time.
\end{thm}
The proof of \cite{SC98} is based on the analysis of the integral equations for the Navier-Stokes equations with the aid of the Cauchy-Kowalewski theorem. On the other hand, quantifying the production of vorticity by the boundary at finite viscosity and its interaction with convective terms is crucial in understanding the behavior of the viscous fluid near impermeable walls. The next result from  \cite{M14} shows that the zero-viscosity limit is verified for a short time by using the vorticity formulation in the half plane, identified with $\RR^2_+$, as long as the initial vorticity stay bounded away from the boundary.
\begin{thm}\label{thm.Maekawa} Let $\omega_0={\rm curl}\, u_0={\rm curl}\, u_0^\nu$ be the initial vorticity for the Euler and Navier-Stokes flows. Assume that
\[
d_0:={\rm dist}\, (\partial \RR^2_+,{\rm supp} \,\omega_0)>0.
\]
Define $u^\nu_{\text{approx}}=u^0 +u^\nu_P$, where $u^\nu_P$ is the boundary layer corrector,
which satisfies a modified Prandtl equation.
Then,
\[
    \|u^\nu - u^\nu_{\text{approx}} \|_{L^\infty((0,T)\times \RR^2)} \leq C\,\nu^{\frac12},
\]
for some constant $C>0$ independent of $\nu$. The time $T$ can be estimated from below as follows:
\[
T\geq c\, \min\{d_0,1\},
\]
for some positive constant $c$ which depends only on $\|\omega_0\|_{W^{4,1}\cap W^{4,2}}$.
\end{thm}
Note that the class of natural test functions $C_{0,\sigma}^\infty (\mathbb{R}^2_+)$ (or even $C_{0,\sigma}^k (\mathbb{R}^2_+)$ for large $k$) is admissible for the  initial data in Theorem \ref{thm.Maekawa}.
This case has been excluded in Theorem \ref{thm.SC98} due to the analyticity
assumption in the entire half space. In Theorem \ref{thm.Maekawa} the $L^\infty$
choice for the norm in which to take the limit is more natural than the energy
norm, as the vorticity-velocity formulation is used to obtain uniform bounds in
viscosity and, as already remarked, the velocity obtained from the Biot-Savart
law, is not in $L^2$ unless the integral of vorticity is zero (see e.g.
\cite{MajdaBertozzi}). Recently, the result of \cite{M14} is extended in
\cite{FTZpreprint} to the three-dimensional case, where the proof is based on
a direct energy method.

The data discussed in Theorems \ref{thm.SC98} and \ref{thm.Maekawa} have
analytic regularity at least near the boundary. To remove the condition of
analyticity is a big challenge,. Recently,  Prandtl expansion has been verified
in \cite{GMMpreprint} in a Gevrey space around a shear flow with a monotone and
concave boundary layer profile. By virtue of the spectral instability shown in
\cite{GGNpreprint'}, the requirement of Gevrey regularity is considered to be
optimal for the stability of the boundary layer at least at a linear level.

If the analyticity is totally absent, one should not expect an asymptotic expansion of the form
\eqref{e.noslipPrandtlExp} to be valid in general, given the underlying strong
instability mechanism at high frequencies. A  classical instability result is
given in \cite{Grenier00}, and recalled in Theorem \ref{thm.Grenier} below,
where the invalidity of \eqref{e.noslipPrandtlExp} is shown when the initial
boundary layer profile is linearly unstable for the Euler equations.

\begin{thm}\label{thm.Grenier} Let $U_s (z)$ be a smooth shear layer satisfying $U_s(0)=0$, such that $U_s {\bf e}_1=(U_s,0)$ is a  linearly unstable stationary solution to the Euler equations.  Let $n$ be an integer, arbitrarily large. Then there exists $\delta_0>0$ such that the following statement holds.
For every large $s$ and sufficiently small $\nu$ there exist $T_\nu>0$ and $v_0^\nu\in H^s (\mathbb{R}^2_+) \cap L^2_\sigma (\mathbb{R}^2_+)$ such that
\begin{align*}
\lim_{\nu\rightarrow 0} T_\nu =0\,, \qquad  \|v_0^\nu \|_{H^s(\mathbb{R}^2_+)}\leq \nu^n\,,
\end{align*}
and the solution $u^\nu$ to \eqref{e.NS} in $\mathbb{R}^2_+$ (under the no-slip boundary condition) with the initial data $u_0^\nu (x) =U_s (\frac{x_2}{\sqrt{\nu}} ) {\bf e}_1 + v_0^\nu (x)$ satisfies the estimate
\begin{align*}
\lim_{\nu\rightarrow 0} \| {\rm curl}\, u^{\nu} (T_\nu) - {\rm curl}\, \big (u_s (T_\nu, \frac{\cdot}{\sqrt{\nu}}) {\bf e}_1 \big )\|_{L^\infty (\mathbb{R}^2_+)} & = \infty\,,\\
\lim_{\nu \rightarrow 0} \| u^\nu (T_\nu) - u_s (T_\nu, \frac{\cdot}{\sqrt{\nu}}) {\bf e}_1 \|_{L^\infty (\mathbb{R}^2_+)} & \geq \delta_0 \nu^\frac14\,.
\end{align*}
Here $u_s (t, z)$ is the smooth solution to the heat equations $\partial_t u_s - \partial_z^2 u_s =0$, $u_s|_{t=0}= U_s$, and $u_s|_{z=0} =0$.
\end{thm}

In the proof given in  \cite{Grenier00}, the time $T_\nu$ is of the order $\mathcal{O} (\nu^\frac12 |\log \nu|)$.
In particular, Theorem \ref{thm.Grenier} implies that the expansion \eqref{e.noslipPrandtlExp} may cease to be valid in general at least in this very short time. Theorem \ref{thm.Grenier} is proved based on the instability of the shear profile $U_s$ for the Euler equations, but in the rescaled variables $X=\frac{x}{\sqrt{\nu}}$.
Recently, spectral instability has been established \cite{GGNpreprint'} even in
the case for which $U_s$ is (neutrally) stable for the Euler equations. This
result indicates the invalidity of Prandtl expansion in the Sobolev framework,
no matter whether the boundary layer profiles possess a good shape, although in
this case there is still no rigorous proof about the nonlinear instability in
Sobolev classes.
On the other hand, as it will be mentioned in Section
\ref{subsubsec.illposed.Prandtl}, the invalidity of the asymptotic estimate
\eqref{e.noslipPrandtlExp} can be also derived from the high-frequency
instability of the shear profile in the Prandtl equations, which is proved in
\cite{GN11}.
Note that, however, the invalidity of \eqref{e.noslipPrandtlExp} observed in
\cite{Grenier00} and \cite{GN11} relies on the assumption that the boundary
layer is formed already at the initial time (and thus, the initial data for the
Navier-Stokes flows is also assumed to depend on the viscosity coefficient, more
precisely, on the fast variable $\frac{x_2}{\sqrt{\nu}}$). It is still open
whether \eqref{e.noslipPrandtlExp} can be disproved, or proved, in the case when
the Sobolev initial data of the Navier-Stokes flows is taken independently in
$\nu$, for in this case there is no boundary layer at the initial time and the
layer  forms only at a  positive time.

Mathematically, the main difficulty in the case of the no-slip boundary condition is the lack of {\em a priori} estimates on strong enough norms to pass to the limit, which in turns is due to the lack of a useful boundary condition for vorticity or pressure. Then, the other types of results found in the literature can be roughly divided into two groups:
\begin{enumerate}[label=(\Roman{enumi}),ref=(\Roman{enumi})]
\item conditional convergence results for generic flows under conditions on the flow that control the  growth of gradients in the layer, such as Kato's condition on the energy dissipation rate, discussed below; \label{item.criteria}
\item convergence results for specific classes of flows, where some conditions as in \ref{item.criteria} are valid automatically, such as parallel flows in pipes and channels discussed below.
\end{enumerate}

Again, for general initial conditions, the zero viscosity limit is sought to hold on the interval of existence of the Euler solutions.

Kato \cite{K84} realized that the vanishing of energy dissipation in a small layer near the boundary is equivalent to the validity of the zero-viscosity limit in the energy space. In fact, this condition is enough to pass to the limit in the non-linear term in the weak formulation of the equations.

\begin{thm}[Kato's criterion]
Let $u^\nu$ be a Leray-Hopf weak solution of the three-dimensional Navier-Stokes
equations \eqref{e.NS} with initial data $u^\nu_0\in L^2(\Omega)$, $\Omega
\subset \mathbb{R}^3$. Let $u^0$ be a strong
solution of the Euler equations \eqref{e.EE} with initial data $u_0\in H^s$,
$s>5/2$ on the time interval $[0,T]$. Assume $u^\nu_0\to u_0$  strongly in
$L^2(\Omega)$.  Then $u^\nu \to u^0$ strongly in
$L^\infty([0,T], L^2(\Omega))$, i.e.,   the vanishing viscosity limit holds on
$,[0,T]$ if and only if, for $T' \leq T$,
\[
     \lim_{\nu \to 0_+} \nu \int_0^{T'} \|\nabla u^\nu(t)\|_{L^2(\Gamma_{c\,\nu})}^2\, {\rm d} t =0,
\]
where $c>0$ is a fixed, but arbitrary, constant and $\Gamma_{c\,\nu}$ is a boundary strip of width $c\,\nu$.
\end{thm}

Variants of Kato's criterion have been established, involving only the partial gradient of the velocity field  and allowing for non-zero boundary velocity \cite{TW1998,W2001}, or involving the vorticity \cite{K07} for instance.
 If the Euler flow satisfies a sign condition at the boundary,
namely, Oleinik's monotonicity condition, that on a half plane reads\
$u^0_1(x_1,0,t)\geq 0$, so that no back flow occurs, then it is enough for the
vorticity to be not too negative in a Kato-type boundary layer (of width $\nu$)
\cite{CKV15}.
Unfortunately, it is not known whether flows generically satisfy Kato's
criterion  at least for short time. In fact,  Kato's criterion cannot hold if
boundary layer separation occurs by a result of  \cite{K08,K14} stated below,
where
it is shown that the zero viscosity limit holds if and only if vorticity
accumulates only at the boundary as a conormal distribution.

\begin{thm}
In the hypothesis of Kato's criterion, the following are equivalent:
\begin{enumerate}[label=(\alph{enumi}),ref=(\alph{enumi})]
\item $u^\nu \to u^0$ in $L^\infty([0,T],L^2(\Omega))$; \label{item.Kell1}
\item $\displaystyle \omega^\nu \to \omega^0 + u^0\times {\bf n}\,
\mu$ weakly-$\ast$ in $L^\infty([0,T),H^{1}(\Omega)')$. \label{item.Kell2}
\end{enumerate}
where  ${\bf n}$ is the unit outer normal to the boundary and $\mu$ is a Radon
measure that agrees with surface area on $\partial\Omega$.
\end{thm}

In two space dimensions, condition \ref{item.Kell2} can be
equivalently restated as:
\[
  \omega^\nu \to \omega^0 - u^0\cdot {\boldsymbol{\tau}}\,
\mu \quad \text{weakly-$\ast$ in } L^\infty([0,T),H^{1}(\Omega)'),
\]
where $\boldsymbol{\tau}$ is the unit tangent vector to $\partial\Omega$
(obtained by rotating {\bf n} counterclockwise by $90^0$ degrees).
It should be noted that the equivalence of \ref{item.Kell1} and \ref{item.Kell2}
above is purely kinematic, in the sense that is a direct consequence of the
validity of the limit and results on weak convergence of gradients. If dynamics
is taken into account, for example when the initial-boundary-value problem for
vorticity is solvable, the convergence in \ref{item.Kell2} can often be improved
to convergence in the sense of Radon measures on the closure of $\Omega$.
In fact, it is enough that the vorticity be uniformly bounded in viscosity
in $L^\infty([0,T),L^1(\Omega))$ (\cite[Corollary 4.1]{K14}). This
will be the case for parallel pipe and channel flows discussed later in this
subsection. In this situation, one can interpret the extra measure on the
boundary appearing in the limit as a {\em vortex sheet} due to a jump in
velocity across the boundary even if there is no fluid outside the domain
$\Omega$.

A consequence of the theorem is that boundary layer separation cannot occur if
the zero-viscosity limit holds, even though the convergence is in a relatively
weak norm, the energy norm, because then the convergence of $u^\nu$ to $u^0$ is
in the $H^1$-Sobolev norm in the interior, and this strong convergence is
incompatible with layer separation.

The vanishing viscosity limit and associated boundary layer can be studied for special classes of flows that satisfy strong symmetry assumptions. In this situation, no additional assumptions are made on the flow, except assuming symmetry of the initial data, as the symmetry is preserved by the Navier-Stokes and Euler evolution, at least for strong solutions.
(For a discussion of possible symmetry breaking in the context of weak solutions, the reader is referred to \cite{BLNNT13}.)
The classes of flows that can be studied are so-called {\em parallel flows} in straight infinite channels or straight infinite circular pipes. These flows can be thought of as generalization of the classical Poiseuille and Couette flows, but they are unsteady and generally non-linear. In fact, the walls of the pipe or channel are allowed to move rigidly along itself, as in the classical Taylor-Couette case, so that the no-slip boundary condition for solutions to \eqref{e.NS} takes the form \eqref{e.noslipBC} with $V\ne 0$.
Parallel channel and pipe flows were considered before in the context of the zero viscosity limit by \cite{W2001}, who lists them as cases for which Kato's criterion  applies.
In fact, it is easy to see that the criterion applies in the extension due to
\cite{TW1998} if the boundary velocity $V$ is not too rough. However, due to the
symmetry in the problem it is possible to obtain a detailed analysis of the
flow in the boundary layer and quantify vorticity production even in the case of
impulsively started and stopped boundary motions, where $V$ is of bounded
variation in time. It should be noted that the boundary layer is not weak here
and, in fact, it has the width predicted by the Prandtl theory proportional to
$\sqrt{\nu}$.
But, because of symmetry, the flow stays laminar and the boundary layer never detaches. A similar analysis for truly non-linear, symmetric flows, such as axisymmetric flows (without swirl) and helical flows seems out of reach at the moment.

In what follows, $\{e_r,e_\phi\}$ will denote the orthonormal frame associated to polar coordinates $(r,\phi)$ in the plane, and  $\{e_r,e_\phi, e_x\}$ will denote the orthonormal frame associated to cylindrical coordinates $(r,\phi,x)$ in space. The symmetric flows that have been considered in the literature are:
\begin{enumerate}[label=(\roman{enumi})]
 \item {\em Circularly Symmetric Flows}: planar flows in a
    disk centered at the origin \, $\Omega=\{ x^2+y^2 < R\}$. The velocity is
of the form
\begin{equation}
   u = V(t) e_\phi,
\end{equation}
 using polar coordinates, where $V(t)$ is a radial function. The vorticity,
which can be identified with a scalar for planar flows, is also radial.

    \smallskip

    \item {\em Plane-parallel flows}: 3D flows in a
    infinite channel, with periodicity imposed  in the $x$ and
    $y$-directions. The velocity takes the form:
        \begin{equation}
            u = (u_1 (t,z), \, u_2 (t,x, z), \, 0),
        \end{equation}
      and is given on the domain
        \begin{align*}
            \Omega := (0, L)^2 \times (0, h)
        \end{align*}
        where  $h$ is the width of the channel and $u_1$, $u_2$ satisfy periodic
       boundary conditions in  $x$ and $y$. The boundary is identified
with the set $\pa \Omega = [0,L]^2 \times [0, h]$

    \smallskip

    \item {\em Parallel pipe flows}:  3D flows in a
      infinite straight, circular pipe, with periodicity imposed along the pipe axis, identified with the $x$-axis.
     The velocity is of the form
        \begin{equation} \label{e.pipeSym}
            u=
                u_{\phi} (t,r) \e_{\phi}
                    + u_x (t,\phi, r) \e_x,
        \end{equation}
      using cylindrical coordinates on the domain
        \begin{align*}
            \Omega := \{(x, y, z) \in \R^3 \, \mid \, y^2 + z^2 < R, \; 0<x<L\},
        \end{align*}
      where $R$ is the radius of circular cross-section and $u_\phi$, $u_x$
    satisfy periodic boundary conditions in $x$.
   The boundary is identified with the set $\pa\Omega =
\{(y,z)\in \RR^2 \,\mid \, y^2+z^2= R\} \times [0, h]$.
\end{enumerate}
For all these flows, the divergence-free condition is automatically satisfied.
In the case of circular symmetry, the Navier-Stokes equations reduces to a heat equation and the Euler flow is steady, making this more of a pedagogical example.
Both for the channel and pipe geometry, symmetry and periodicity ensure
uniqueness of solutions to the  Navier-Stokes equations and the Euler
equations, in particular by forcing the only
pressure-driven flow to be the trivial flow.
The well-posedness is global in time for both Euler and Navier-Stokes for sufficiently regular initial data.

As an illustration, only parallel pipe flows will be discussed here, which is the most interesting case due to the effect of curvature of the boundary. The reader is referred to \cite{Matsui1994,BW02,LMN08,LMNT08} for the case of circularly symmetric flows, and to \cite{MT08,MNW10} for the case of channel flows.
In parallel pipe flows, the velocity is independent of the variable along the
pipe axis and, in any circular cross section of the pipe, it is the sum of a
circularly symmetric, planar velocity field and a velocity field pointing in the
direction of the axis. As in the case of plane-parallel flows, even though the
flows are not planar,
the Navier-Stokes and Euler equations reduce
to a weakly non-linear system in only two space variables, given respectively by
(for simplicity it is assumed that the boundary is stationary):
\[
   \begin{cases}
                               \displaystyle        \frac{\pa u^\nu_{\phi}}{\pa t}
                                        - \nu \Delta u^\nu_{\phi}
                                        + \nu \frac{1}{r^2} u^\nu_{\phi}
                                        = 0
                                        &
                                        \text{in } (0,T)\times \Omega\,,\\

                                       \displaystyle  \frac{\pa u^\nu_{x}}{\pa t}
                                        - \nu \Delta u^\nu_x
                                        + \frac{1}{r} u^\nu_{\phi} \frac{\pa u^\nu_x}{\pa \phi}
                                        = 0
                                        &
                                        \text{in } (0,T) \times \Omega\,,\\
                                        \displaystyle - \frac{1}{r} (u^\nu_{\phi})^2
                                        + \frac{\pa p^\nu}{\pa r}
                                        = 0
                                    &
                                       \text{in } (0,T) \times \Omega\,,\\
                                       u^\nu_i
                                       = 0\,, \text{ } i=\phi, x
                                       &
                                       \text{on } (0,T) \times \pa\Omega\,,\\
        \end{cases}
  \]
and by
\[
   \begin{cases}
                                        \displaystyle \frac{\pa u^0_{\phi}}{\pa t}=0
                                          &
                                        \text{in } (0,T) \times \Omega\,,\\

                                        \displaystyle \frac{\pa u^0_{x}}{\pa t}
                                        + \frac{1}{r} u^0_{\phi} \frac{\pa u^0_x}{\pa \phi}
                                        = 0
                                        &
                                        \text{in } (0,T) \times \Omega\,,\\
                                        \displaystyle - \frac{1}{r} (u^0_{\phi})^2
                                        + \frac{\pa p^0}{\pa r}
                                        = 0
                                        &
                                        \text{in } (0,T) \times \Omega\,.
        \end{cases}
  \]
The Navier-Stokes system is amenable to the analysis of the vanishing viscosity limit primarily because it is diffusion dominated and because the pressure is slaved to the velocity and drops out of the momentum equation.

A detailed analysis of the boundary layer using techniques borrowed from semiclassical analysis was performed in \cite{MT11}, for ill-prepared data. There, in particular, convergence rates in viscosity for the $L^\infty$ norm were derived by constructing a  parametrix to a suitable associated linear problem and taking the corrector as the double layer potential associated to this problem. For well-prepared data, convergence rates in higher Sobolev norms were obtained by the use of flow correctors and effective equations in \cite{HMNW11}. By the use of a different type of correctors, it is possible to obtain similar results for ill prepared  data and quantify production of vorticity at the boundary \cite{GKLMN15}.
Similarly to the case of Navier boundary conditions, one defines an approximate Navier-Stokes solution $u^\nu_\text{approx}$ as a sum of an outer solution $u^\nu_{\text{ou}}$ and an inner solution $u^\nu_{\text{in}}$ . At zero order in viscosity, $u^\nu_{\text{ou}}=u^0$, the Euler solution, while $u^\nu_{\text{in}}$ is given by a smooth, radial  cut-off $\psi$ supported in a collar neighborhood of the boundary times a corrector $\theta$ of the form:
\begin{equation}\label{e.pipeCorrector}
   \theta(t,x) = \theta_{\phi} (t,r) \e_{\phi}
                    +
                    \theta_{x} (t,\phi, r) \e_x\,,
\end{equation}
using again cylindrical coordinates $(r,\phi,x)$, where  $\theta_{\phi}$ and $\theta_x$ satisfy weakly coupled parabolic systems. Then, the following convergence rates can be obtained \cite{MT11,HMNW11,GKLMN15}.

\begin{thm}
Assume $u_0\in H^k(\Omega)$, $k$ large enough ($k>4$ suffices), and has symmetry \eqref{e.pipeSym}. Then, the zero-viscosity limit hold on $(0,T)$, for all $0<T<\infty$ and, in particular:
    \begin{equation}\label{e.pipeConv}
        \left\{
                \begin{array}{l}
                        \| u^\nu- u^\nu_\text{approx}\|_{L^{\infty}(0, T; L^2(\Omega))}
                        + \nu^{\frac{1}{2}} \| \nabla u^\nu- \nabla u^\nu_{\text{approx}} \|_{L^{2}(0, T; L^2(\Omega))}
                                \leq \kappa_T \, \nu^{\frac{3}{4}},\\
                              \| u^\nu - u^0 \|_{L^{\infty}(0, T;
L^2(\Omega))} \leq \kappa_T \, \nu^{\frac{1}{2}},
                \end{array}
        \right.
        \end{equation}
where $\kappa_T$ is a constant independent of $\nu$.
In addition,
\begin{align}\label{e.pipeVortConv}
    \omega^\nu \to \omega^0 + (u^0\times {\bf n}) \mu
        \quad \text{weakly}^* \text{ in } L^\infty(0, T; \mathcal{M}(\Bar{\Omega})),
\end{align}
where $\mathcal{M}(\Bar{\Omega})$ is the space of Radon measures on $\Bar{\Omega}$, and $\mu$ is a measure supported on $\pa\Omega$, which on the boundary agrees with the normalized surface area.
\end{thm}

A complication over plane-parallel flows is that the effect of
non-vanishing curvature cannot be neglected in the analysis.
Furthermore, in cylindrical coordinates
the behavior of the solution near the axis cannot be controlled as well as it can be away from the axis, similarly to the case of axisymmetric flows. To overcome this difficulty,   a
two-step localization, one near the boundary where curvilinear coordinates are
used, the other near the axis where Cartesian coordinates and energy estimates
are employed, is utilized. As a consequence, however, the error estimates for the approximate solution suffer from the loss of one derivative. In particular, the estimates for the correctors are not as sharp as in the case of a pipe with annular cross section.

As seen below, there are other situations for which it is possible to pass to the limit  due to the fact that the boundary layer is weak or absent.
This is the case, for example, of flows outside shrinking obstacles, if the obstacle is shrinking  faster than viscosity vanishes.
For such flows the local Reynolds number, built by taking the size of the obstacle as characteristic length, stays of order one, as already observed in \cite{ILN09}. In this context, the Navier-Stokes solution in the exterior of the obstacle is expected to converge to the Euler solution in the whole space.
Most of the results concern flows in the plane, as the vorticity-velocity formulation is used to obtain the inviscid solution. At the same time, the fact that the exterior of compact obstacles is not simply connected in two space dimensions adds some technical difficulties, which are overcome by assuming that the circulation around the obstacles is zero.

Let $\epsilon$ be the scale of the obstacle. The vanishing viscosity limit was shown to hold in the exterior of one obstacle diametrically shrinking to a point in \cite{ILN09} by assuming the condition \, $\epsilon\leq C\,\nu$ for some positive constant $C$, which depends on the initial data for the Euler equations in $\RR^2$, $u_0$, and the shape of the obstacle, and assuming that the initial condition for the Navier-Stokes solution, $u^\nu_0$, extended by zero to the whole plane, converges to $u_0$ in $L^2(\RR^2)$, with an optimal rate of convergence of $\sqrt{\nu}$. (See \cite{KLN09} for the opposite situation of an expanding domain.)
This result can be extended to the exterior of a finite  number of fixed
obstacles. It is interesting to ask whether a similar result hold in the setting
of a porous medium,  that is, if the domain for the viscous flow is the exterior
of an array of particles. It is known that homogenization of the Navier-Stokes
equations and Euler equations gives different filtration laws (Darcy or
Brinkman, for example) depending on the relative ratio the particle size
$\epsilon$ and inter-particle distance $d_\epsilon$, and the permeability of the
homogenized medium is very different between the viscous and inviscid cases (see
the discussion in \cite{LCM15,MikelicPaoli} and references therein). Therefore,
it is relevant to study the joint limit of vanishing $\epsilon$, $d_\epsilon$,
$\nu$.

In \cite{LCM15}, the limit was established under the condition that $d_\epsilon
> \epsilon$ and $\epsilon\leq A\,\nu$; see Theorem \ref{thm.LCM15} below. In
this regime, one expects that the limit Euler flow, defined in the whole plane,
does not feel the presence of the porous medium. Below, for each $\epsilon$, the
domain $\Omega_\epsilon$ is set as the viscous fluid domain, which, for
simplicity, is  the exterior of a regular array of identical particles arranged
in a square.
\begin{thm}\label{thm.LCM15}
Given $\omega_{0}\in C^\infty_{c}(\RR^2)$,  let $u^0$ the solution of the Euler
equations  in the whole plane with initial condition $u_{0}:=K_{\RR^2}[\omega_0] $. For any $\epsilon,\nu>0$, let $d^\epsilon
\geq \epsilon$ and let  $u^{\nu,\epsilon}$ be the solution of  the
Navier-Stokes equations in
$\Omega^\epsilon$ with initial velocity $u_{0}^{\nu,\epsilon}$. Then, there
exists a constant $A$ depending only on  the particle shape, such that if
$$\displaystyle \frac{\epsilon}{d_{\epsilon}} \leq \frac{A \nu}{ \| \omega_{0}\|_{L^1\cap L^\infty(\R^2)}},$$
and if $\omega_{0}$ is supported in $\Omega^\epsilon$, then for any $T>0$ we have
\begin{equation} \label{ineq.main}
         \sup_{0\leq t\leq T}
  \|u^{\nu,\epsilon}-u^0\|_{L^2(\Omega^\epsilon)}\leq B_{T} \left(
\frac{\sqrt\nu}{d_{\epsilon}} +\|u_{0}^{\nu,\epsilon}-u_0\|_{L^2(\Omega^\epsilon)}\right)
\end{equation}
where $B_{T}$ is a constant depending only on $T$, $\|\omega_{0}\|_{L^1\cap
W^{1,\infty}(\R^2)}$,  and the particle shape.
\end{thm}

It is then possible to construct initial data $u^{\nu,\epsilon}_0$ such that
$u^{\nu,\epsilon}_0\to u_0$ in $L^2$, which establishes the limit with rate
$\sqrt\nu/d_{\epsilon}$. Therefore, there is a ghost of the porous medium in the convergence rate.
It should be noted that in the case of the Darcy-Brinkman system, the equations for  the boundary corrector are linear, and thus, the passage to the zero-viscosity limit is  possible \cite{KTW11,HW14}.


\subsection{Non-characteristic boundary case}\label{s.nonchar}

One of the main difficulties in treating the vanishing viscosity limit for classical no-slip boundary conditions is the fact that the boundary is characteristic for the problem, that is it consists of streamlines for both the viscous and inviscid flows. Hence, any attempt to control the flow in the interior from the boundary seems unsuccessful unless analyticity or monotonicity of the  data is imposed.

If non-characteristic boundary conditions are imposed, in particular, if the walls are permeable, then under certain conditions the boundary layer is stable, hence there is no layer separation and one expects the vanishing viscosity limit to hold.
This is the case when injection and suction rates are imposed at the boundary.
For simplicity we describe the set up in the geometry of a (periodized) channel
$[0,L]^2\times [0,h]$, where injection and suction is imposed at the top
and bottom walls.  The velocity at the boundary for Navier-Stokes is then given
as:
\begin{equation} \label{e.noncharBC}
  u^\nu_i (t,x)= (0,0, -U_i (t,x_1,x_2)), \qquad i={\rm top},{\rm bot},
\end{equation}
where $U_i \geq a_i>0$ for some constants $a_i$ and ${\rm top}$, ${\rm bot}$ refer to top and bottom of the channel. For Euler, one needs to specify the entire velocity and the inlet/outlet. These conditions are also appropriate when a domain is truncated, e.g. for computational reasons, when making a Galilean coordinate transformation.

By correcting the velocity field, it was shown in \cite{TWsuction} that the zero viscosity limit holds with sharp rates of convergence of $\nu^{1/4}$ in the uniform norm. In particular, there is only a stable boundary layer at the suction wall (the bottom) that is exponentially small. It is interesting to note that, differently than in the non-linear case, for the Oseen equations, adding injection and suction at the boundary does not seem to change the size of the boundary layer  (see \cite{LSOseen01}).

\subsection{Prandtl equations for the boundary layer}\label{s.Prandtl}

This subsection is devoted to an  overview of the study of the Prandtl equations,
introduced by Prandtl in 1904 in order to describe a viscous incompressible flow near the boundary
at high Reynolds numbers \cite{Pra1904}.
The Prandtl equations are derived from the Navier-Stokes equations with no-slip boundary condition,
and their derivation is briefly recalled here in the case that the fluid domain
is the half plane $\mathbb{R}^2_+$. The initial-boundary-value problem for  the
Navier-Stokes equations reads:
\begin{equation}\label{NS.sec3.4}
  \left\{
\begin{aligned}
 \partial_t u^\nu + u^\nu \cdot \nabla u^\nu + \nabla p^\nu & = \nu \Delta u^\nu \,, \qquad t>0\,, ~~x\in \mathbb{R}^2_+\,,\\
{\rm div}\, u^\nu & =0\,,  \qquad \qquad t\geq 0\,, ~~x\in \mathbb{R}_+^2\,,\\
u^\nu|_{t=0} & = u_0^\nu\,,\qquad \qquad \qquad ~~ x\in \mathbb{R}^2_+ \,, \\
u^\nu|_{x_2=0} & = 0\,,\qquad \qquad \qquad ~~ t>0, \; x_1\in \mathbb{R}.
\end{aligned}\right.
\end{equation}
As is discussed in the previous sections, the equations in the limit $\nu =0$ are
the Euler equations, for which only the impermeability condition $u_2^\nu=0$ on $\partial\mathbb{R}^2_+$ can be prescribed. Heuristically, such an incompatibility in the boundary condition leads to a fast change of the tangential component of the velocity field, which is $u_1^\nu$ when the fluid domain is $\mathbb{R}^2_+$.
As a result, the derivative of $u_1^\nu$ in the vertical direction tends to have a singularity near the boundary and forms a boundary layer.
To study the formation of the boundary layer,  Prandtl made the ansatz that $u^\nu$ near the boundary has the following asymptotic form:
\begin{align}\label{eq.ansatz}
u_1^\nu (t,x_1,x_2) \sim u_1^P (t,x_1,\frac{x_2}{\sqrt{\nu}})\,, \qquad u_2^\nu (t,x_1,x_2)  \sim \sqrt{\nu} u_2^P (t,x_1,\frac{x_2}{\sqrt{\nu}})\,.
\end{align}
The thickness of the boundary layer $\mathcal{O}(\sqrt{\nu})$ is coherent with the parabolic nature of the Navier-Stokes equations. The underlying assumption here is that the velocity $u^\nu$ remains of  order $\mathcal{O}(1)$ in all $\partial_1^k u^\nu$, $k=0,1,\cdots$, in the limit $\nu\rightarrow 0$, and then the vertical component $u_2^\nu$ is expected to be of the order $\mathcal{O}(\sqrt{\nu})$ since the boundary condition for the normal component is preserved in the limit, and it is also compatible with the divergence free condition.
By formally substituting the ansatz \eqref{eq.ansatz} into the first equation of \eqref{NS.sec3.4}, the velocity profile $u_1^P$ and the associated pressure $p^P$ should obey the equations
\begin{align*}
\partial_t u_1^P + u^P \cdot \nabla u_1^P + \partial_1 p^P - \partial_2^2 u_1^P  = 0\,, \qquad \partial_2 p^P=0\,,
\end{align*}
and $u_1^P$ must satisfy the no-slip boundary condition $u_1^P=0$ on $\partial\mathbb{R}^2_+$.
Here the spatial derivatives are for the rescaled variables $X_1=x_1$ and $X_2=\frac{x_2}{\sqrt{\nu}}$,
but in this subsection $\partial_j$ will denote both
$\frac{\partial}{\partial x_j}$ and $\partial_j=\frac{\partial}{\partial
X_j}$ for notational ease. The rescaled variables $X$ will be also
relabeled as $x$ from now on.
The vertical component $u_2^P$ is recovered from $u_1^P$ and the boundary
condition $u_2^P=0$ on $\partial\mathbb{R}^2_+$ by virtue of the second
equation (divergence-free condition) in \eqref{NS.sec3.4}, which yields
\begin{align*}
u_2^P  (t,x) = -\int_0^{x_2} \partial_{1} u_1^P (t,x_1,y_2) \dd y_2\,.
\end{align*}
The velocity  in the boundary layer has to match with the outer flow which is
assumed to satisfy the Euler equations. This requirement leads to
the following boundary condition  ({\em matching conditions}) on $u_1^P$ and
$p^P$ at $x_2=\infty$:
\begin{align*}
\lim_{x_2\rightarrow \infty} u_1^P = u^E\,, \qquad \lim_{x_2\rightarrow \infty} p^P = p^E\,,
\end{align*}
where $u^E (t,x_1) = u^0(t,x_1,0)$ and $p^E (t,x_1)= p^0 (t,x_1,0)$, and $(u^0,p^0)$ is the solution to the Euler equations \eqref{e.EE} in $\mathbb{R}^2_+$. Since $p^P$ must be independent of $x_2$, because of the  equation $\partial_2p^P=0$ in $\mathbb{R}^2_+$, the matching condition on $p^P$ at $x_2=\infty$ implies that
\begin{align*}
p^P = p^E\,.
\end{align*}
That is, the pressure field is not one of the unknowns  in the Prandtl equations.
Collecting the above equations gives the Prandtl equations in $\mathbb{R}^2_+$ (in the spatial variables),
which are a system of equations  for the scalar unknown function $u_1^P$:
\begin{equation}\label{Pra.sec3.4}
  \left\{
\begin{aligned}
& \partial_t u_1^P + u^P \cdot \nabla u_1^P  - \partial_2^2 u_1^P  = - \partial_1 p^E \,, \\
& u_2^P  = -\int_0^{x_2} \partial_{1} u_1^P  \dd y_2\,, \\
& u_1^P|_{t=0} = u^P_{0,1}\,, \qquad u_1^P|_{x_2=0} =0\,, \qquad  \lim_{x_2\rightarrow \infty} u_1^P = u^E\,. \\
\end{aligned}\right.
\end{equation}
The reader is referred to \cite{Schlichting,Mas2007,LiWa16} for more details
about the formal derivation of the Prandtl equations. Note that, by
taking the boundary trace in the Euler equations, the data $(u^E, p^E)$ coming
from the Euler flows in the outer region is subject to the Bernoulli law
\begin{align}\label{eq.Bernoulli.sec3.4}
\partial_t u^E + u^E \partial_1 u^E + \partial_1 p^E = 0\,.
\end{align}
The Prandtl equations are deceptively simpler than
the original Navier-Stokes equations. In fact , due to the inherent instability of boundary layers, well-posedness of \eqref{Pra.sec3.4} has been proven only in some specific situations (Section \ref{subsubsec.wellposed.Prandtl}), while strong ill-posedness results are present in the literature (Section \ref{subsubsec.illposed.Prandtl}).

\subsubsection{Well-posedness results for the Prandtl equations}\label{subsubsec.wellposed.Prandtl}

The Prandtl equations are known to be well-posed under some restricted conditions.
This subsection gives a list of the categories in which the well-posedness of the Prandtl equations  holds, at least for short time.

\vspace{0.3cm}

\noindent
{\bf (I) Monotonic data.} This category is the most classical in the theory of the Prandtl equations. The system \eqref{Pra.sec3.4} is  considered for $0<t<T$ and for $(x_1,x_2)\in \Omega_1\times \mathbb{R}_+$,
where $\Omega_1$ is usually set as either $\{0<x_1<L\}$, $\mathbb{T}$, or $\mathbb{R}$. When $\Omega_1=\{0<x_1<L\}$, an additional boundary condition has to be imposed on $u_1^P$ at the boundary $\{x_1=0\}$:
\begin{align*}
u_1^P (t,0,x_2) = u_{1,1}^P (t,x_2)\,.
\end{align*}
The given boundary data $u_{1,1}^P$ also has to be compatible with the monotonicity.
The basic assumption describing the monotonicity is
\begin{align}
& \partial_2 u_{0,1}^P (x_1, x_2)>0\,, \qquad x_1\in \overline{\Omega_1}\,, ~ x_2 \geq 0\,, \label{monotone.1}\\
& \partial_2 u_{1,1}^P (t,x_2) >0\,, \qquad ~~ t>0\,, ~ x_2\geq 0\,.\label{monotone.2}
\end{align}
As a compatibility condition, the outer flow $u^E$ and given data $u_{0,1}^P$, $u_{1,1}^P$ must be positive for $x_2>0$, and the solution $u_1^P$ is also expected to be positive for $x_2>0$ together with its derivative in the $x_2$ variable.
The solvability of \eqref{Pra.sec3.4} in this class has been established by Oleinik and her co-workers, especially for the case $\Omega_1=\{0<x_1<L\}$.( See \cite{Oleinik63, Oleinik66, Oleinik66Stability}. The reader is also referred to \cite{OleinikSamokhin99} for more details and references.)
The steady problem is solved in \cite{Oleinik63} for small $L>0$,
and this local existence result is extended in \cite{MaShi1984},
where it is shown that the solution can be continued to the separation point.
For the unsteady problem, unique solvability is proved in \cite{Oleinik66} for a short time if $L$ is given and fixed, while for an arbitrary time if $L$ is sufficiently small.
The stability of the steady solutions is shown in \cite{Oleinik66Stability}.

A natural question arises, already present in the monograph \cite{OleinikSamokhin99}, namely,
under which condition the solutions exist globally in time without any smallness of $L>0$.
A significant contribution to this problem is given by the work \cite{XZ04},
where the global existence of weak solutions to \eqref{Pra.sec3.4} is proved
when the pressure gradient is favorable:
\begin{align}\label{pressure.condition}
\partial_1p^E (t,x_1) \leq 0\,, \qquad t>0\,, ~~ 0<x_1<L\,.
\end{align}
The analysis in \cite{Oleinik66,XZ04} uses the classical Crocco transformation
\begin{align*}
\tau = t\,, \qquad \xi = x_1\,, \qquad \eta = \frac{u_1^P(t,x_1,x_2)}{u^E(t,x_1)}\,, \qquad w (\tau,\xi,\eta) = \frac{\partial_2u_1^P (t,x_1,x_2)}{u^E (t,x_1)}\,,
\end{align*}
which transforms the domain $\{(t,x_1,x_2)~|~0<t<T\,, ~0<x_1<L\,,~x_2>0\}$, $T>0$, into
\begin{align*}
Q_T = \{ (\tau,\xi,\eta)~|~0<\tau<T\,, ~0<\xi<L\,, ~0<\eta<1\}\,.
\end{align*}
Then the Prandtl equations for the case $\Omega_1 = \{0<x_1<L\}$ are transformed
into the system
\begin{equation}\label{CPra.sec3.4}
  \left\{
\begin{aligned}
& \partial_\tau w^{-1} +\eta u^E \partial_\xi w^{-1} + A \partial_\eta w^{-1} - B w^{-1}    = -\partial_\eta^2 w \qquad {\rm in}~~Q_T \,, \\
& w|_{\tau=0} =w_0 = \frac{\partial_2u_{1}^P}{u^E}|_{t=0}\,, \quad w|_{\xi=0}=w_1\,, \quad (w\partial_\eta w)|_{\eta=0}= \frac{\partial_1 p^E}{u^E}\,, \quad w|_{\eta=1} = 0\,.
\end{aligned}\right.,
\end{equation}
where
\begin{align*}
A = (1-\eta^2)\partial_1 u^E + (1-\eta)\frac{\partial_t u^E}{u^E} \,, \qquad B= \eta \partial_1 u^E + \frac{\partial_t u^E}{u^E} \,, \qquad w_1 = \frac{\partial_2 u_{1,1}^P}{u^E}|_{x_1=0}\,.
\end{align*}
The following result is proved in \cite[Theorem 1.1]{XZ04}.

\begin{thm}\label{thm.monotone} Assume that \eqref{monotone.1} and \eqref{monotone.2} hold together with compatibility conditions.
If, in addition, the pressure condition \eqref{pressure.condition} holds,
then there exists a weak solution $w\in BV(Q_T)\cap L^\infty (Q_T)$ to
\eqref{CPra.sec3.4} such that for some $C>0$,
\begin{align*}
C^{-1} (1-\eta) \leq w \leq C (1-\eta) \qquad {\rm in}~~Q_T\,,
\end{align*}
and $\partial_\eta^2w$ is a locally bounded measure in $Q_T$.
\end{thm}

In Theorem \ref{thm.monotone}, the initial and boundary conditions are satisfied
in the sense of the trace,
and \eqref{CPra.sec3.4} is considered in the sense of distributions. The reader
is referred to \cite{XZ04} for more on the properties and the regularity of weak
solutions obtained in this theorem.
As a consequence of Theorem \ref{thm.monotone}, global existence of weak solutions to \eqref{Pra.sec3.4} follows. However, uniqueness and smoothness of  weak solutions in the Crocco variables seem to be still unsettled. In particular, when $L$ is not small enough the global existence of smooth solutions to the Prandtl equations remains open even under the monotonicity conditions \eqref{monotone.1}, \eqref{monotone.2}, and the pressure condition \eqref{pressure.condition}.

The Crocco transformation has been a basic tool in the classical works \cite{Oleinik66,Oleinik66Stability,OleinikSamokhin99,XZ04}.
Recently, an alternative, independent, approach has been presented in
\cite{AWXY2015,MaWo2015},
where the crucial part of the proof is based on a direct energy method but for
new dependent variables.
In particular, the Crocco transformation in not needed in this new approach. The key new unknown is
$w=\displaystyle \partial_2 (\frac{u_1^P}{\tilde \omega^P})$ in \cite{AWXY2015} and $g=\displaystyle \omega^P \partial_2 (\frac{u_{1}^P}{\omega^P})$ in \cite{MaWo2015}, where $\omega^P=\partial_2 u_1^P$ and $\tilde \omega^P=\partial_2 \tilde u_1^P$ represent the vorticity fields of the Prandtl flow $u_1^P$ and of the background Prandtl flow $\tilde u_1^P$, respectively. The function $w$ is introduced in \cite{AWXY2015} in the analysis of the linearized Prandtl equations around $\tilde u_1^P$, and the function $g$ is analyzed in \cite{MaWo2015} for the nonlinear energy estimate. As explained in \cite{MaWo2015} for example,  the crucial obstacle one meets in the energy estimate  for the standard unknowns $\partial_1^j u_1^P$ or $\partial_1^j \omega^P$ is the presence of the terms $(\partial_1^j u_2^P) \omega^P$ or $(\partial_1^j u_2^P) \partial_2 \omega^P$, since they contain the highest order derivatives in $x_1$ and do not vanish after integration by parts. The new unknown is in fact chosen so that these crucial terms cancel.
The development of the Prandtl theory in \cite{AWXY2015,MaWo2015} has had  a significant impact in the field,
and has led to significant recent progress  \cite{KMVW14,D-VM15,IV2016}.

Very recently, the Prandtl equations have been studied in the three-dimensional
half space in  \cite{LWY16,LWYPreprint}. These works show that the monotonicity
condition on the tangential velocities is not enough to ensure  local
well-posedness, and a sharp borderline condition is found for
well-posedness/ill-posedness.

\vspace{0.3cm}

\noindent
{\bf (II) Analytic data.} The second classical category for local
well-posedness of the Prandtl equations is the space of analytic functions.
Under analyticity of the initial data and of the outer Euler flow, local
existence of the Prandtl equations has been  proven in \cite{SC98} after the
pioneer work \cite{Asano88,Asano88'}, where analyticity is imposed on both
variables $x_1$ and $x_2$. Later, it was realized that the condition of
analyticity in the vertical variable $x_2$ can be removed, and the local
well-posedness is known to hold under only analyticity in the tangential
variable \cite{LCS03,KV13,CLM13,CLS14}.

To be precise, a typical existence result available by now in this category is stated here. Let $\sigma\in \mathbb{R}$, $\alpha,\beta,T>0$. The space $\mathcal{H}^{\sigma,\alpha}$ is the space of functions $f(x_1,x_2)$, $2\pi$ periodic in $x_1$, such that the norm
\begin{align*}
|f|_{\sigma,\alpha} = \sum_{j\leq 2} \sup_{x_2\in \mathbb{R}_+} \langle x_2\rangle^\alpha \sum_{k\in \mathbb{Z}} |\partial_2^j \hat{f} (k, x_2)|e^{|k|\sigma}\,, \qquad \langle x_2\rangle = (1+x_2^2)^\frac12\,,
\end{align*}
is finite. Here $\hat{f}(k,x_2)$ is the $k$th Fourier mode of $f$ with respect to $x_1$.
The space $\mathcal{H}^{\sigma,\alpha}_{\beta,T}$ is the space of functions $f(t,x_1,x_2)$, $2\pi$ periodic in $x_1$, such that the norm
\begin{align*}
|f|_{\sigma,\alpha,\beta,T} = \sum_{j\leq 2} \sup_{0\leq t\leq T} | \partial^j_2 f(t) |_{\sigma-\beta t,\alpha} + \sup_{0\leq t\leq T}  | \partial_t f (t) |_{\sigma-\beta t,\alpha}
\end{align*}
is finite. The space $\mathcal{H}^\sigma_{\beta,T}$ is the space of functions $f(t,x_1)$, $2\pi$ periodic in $x_1$, such that the norm
\begin{align*}
|f|_{\sigma,\beta,T} = \sum_{i=0,1} \sup_{0\leq t\leq T} | \partial_t^i f(t)|_{\sigma-\beta t}
\end{align*}
is finite. The spaces $\mathcal{H}^{\sigma,\alpha}$ and $\mathcal{H}^{\sigma,\alpha}_{\beta,T}$ are used for the boundary layer profiles, while $\mathcal{H}^\sigma_{\beta,T}$ is used for the Euler flows.
The next theorem is proved in \cite{CLS14}, where local solvability is obtained even with incompatible initial data, that is,  $u_{0,1}^P|_{x_2=0}\ne 0$.

\begin{thm}\label{thm.analyticity} Let $u^E\in \mathcal{H}^{\sigma_0}_{\beta_0,T_0}$ and $u_{0,1}^P-u^E|_{t=0}\in \mathcal{H}^{\sigma_0,\alpha}$ for some $\sigma_0,\beta_0,T_0>0$ and $\alpha>\frac12$.
Then there exist $\sigma\in (0,\sigma_0)$, $\beta\in (0,\beta_0)$, and $T\in (0,T_0)$ such that Prandtl equations
\eqref{Pra.sec3.4} admit a unique solution $u_{1}^P$ in $[0,T]$ of the form
\begin{align*}
u_{1}^P (t,x_1,x_2) =  - 2 u_{0,1}^P (x_1,0) \, {\rm erfc}\, \big (\frac{x_2}{2\sqrt{t}} \big )  + \tilde u (t,x_1,x_2) + u^E (t,x_1)\,,
\end{align*}
where $\tilde u \in \mathcal{H}^{\sigma,\alpha}_{\beta,T}$.
Here $\displaystyle {\rm erfc}\, \big (\frac{x_2}{2\sqrt{t}} \big )= \displaystyle \frac{1}{\sqrt{\pi t}} \int_{x_2}^\infty \exp  \big ( {-\frac{y_2^2}{4t}} \big ) \dd y_2$.
\end{thm}

As in \cite{Asano88,SC98,LCS03,CLM13}, the proof given in  \cite{CLS14} relies on, the abstract Cauchy-Kowalewski theorem, which is applied to the integral equations associated with \eqref{Pra.sec3.4}.
The existence of solutions under a polynomial decay condition on $u_{0,1}^P-u^E|_{t=0}$ has been shown first in \cite{KV13}. The proof of \cite{KV13} is based on the direct energy method, rather than the use of the abstract Cauchy-Kowalewski theorem and integral equations.

Without monotonicity, one cannot expect global existence of smooth solutions of the Prandtl equations to hold in general even if the given data is analytic.  Indeed, in this case, the existence of finite-time blowup solutions can be shown  (\cite{EE97}; see Section \ref{subsubsec.illposed.Prandtl} below).
However, the class of initial data for blowup solutions in \cite{EE97} must have $\mathcal{O}(1)$ size,
and, hence, it is still not clear whether global existence holds  for sufficiently small data or not.
Recently an important progress has been achieved in this direction, and long-time well posedness is established in \cite{ZZ2016,IV2016} for small solutions in the analytic functional framework. In \cite{ZZ2016} the life span of local solutions is estimated from below to be of order $\mathcal{O}(\epsilon^{-\frac34})$ when the uniform outer Euler flow $u^E=\underline{u}$ is of the order $\mathcal{O}(\epsilon^\frac53)$ and the initial data $u_{0,1}^P$ is $\mathcal{O}(\epsilon)$ in a suitable norm measuring  tangential analyticity.
In \cite{IV2016}, almost global existence is established if  the smallness condition on the uniform Euler flow is removed,  and the life span is then estimated from below as $\mathcal{O}(\exp (-\frac{1}{\epsilon \log \epsilon}))$, $0<\epsilon\ll 1$, when $u_{0,1}^P$ is $\mathcal{O}(\epsilon)$.

\vspace{0.3cm}

\noindent
{\bf (III) Other categories: beyond analyticity or monotonicity.} Without monotonicity or analyticity of the given data,
the solvability of the unsteady Prandtl equations becomes a highly difficult problem even locally in time.
There are a few classes of initial data that are not strictly included in the categories (I) and (II) above, but for which the Prandtl equations can be solved for a short time.

\noindent
(1) {\it Gevrey class with a non-degenerate vorticity.} As it will be seen in Section \ref{subsubsec.illposed.Prandtl}, the Prandtl equations are ill posed in general in Sobolev spaces. This ill-posedness is due to the instability at high frequencies for the tangential components, occurring when the monotonicity of the given data is absent. A key argument  for this instability is given in \cite{G-VD10}, where it is proved that the linearization around the non-monotonic shear flow satisfying \eqref{condition.critical} for some $a>0$ has a solution growing exponentially in time with growth rate $\mathcal{O}(|n|^\frac12)$ for high tangential frequencies $n$. Although such a high-frequency instability yields the ill posedness in Sobolev spaces, there is still hope to obtain the well posedness for initial data whose $n$th Fourier mode in the $x_1$ variable decays in  order $\mathcal{O}(e^{-c|n|^\gamma})$ for $|n|\gg 1$ with some $\gamma >\frac12$, that is, the Gevrey class less than $2$. A crucial difference between the Gevrey class $1$ ($\gamma=1$, analytic functions) and the Gevrey classes $m$ with $m>1$ ($\gamma=\frac1m$) is that the latter class can contain compactly supported functions.
The verification of the instability in \cite{G-VD10} motivates the work of \cite{D-VM15}, where the local solvability is established for a set of initial data without monotonicity, but belonging to the Gevrey class $\frac74$ in the $x_1$ variable.
The key condition for the initial data in \cite{D-VM15} is that the monotonicity is absent only on a single smooth curve but in a non-degenerate manner. More precisely, in \cite{D-VM15} it is assumed that $u^E=p^E=0$ and $u_{0,1}^P$ is periodic in $x_1$ with Gevrey $\frac74$ regularity, and that
\begin{align}\label{condition.GVM}
\partial_2 u_{0,1}^P (x_1,x_2) = 0 \quad {\rm iff} ~~x_2 = a_0(x_1)>0 \quad {\rm with} ~~ \partial_2^2u_{0,1}^P  (x_1, a_0(x_1))>0 \qquad {\rm for ~all}~~ x_1\in \mathbb{T}\,.
\end{align}
Note that the condition \eqref{condition.GVM} is a natural counterpart of \eqref{condition.critical}.
The crucial observation for the proof in \cite{D-VM15} is that in the region where the monotonicity is absent the flow is expected to be convex by virtue of the non-degenerate condition, while away from the curve of the critical points one can use the monotonicity of the flow. However, due to the nonlocal nature of the Prandtl equations,
taking advantage of each in the different regions requires an intricate analysis, and this difficulty is overcome in \cite{D-VM15} by introducing various kinds of energy.

\vspace{0.1cm}

\noindent
(2) {\it Data with multiple monotonicity/analyticity regions.}
The flows in this class are introduced in \cite{KMVW14},
where the local existence and uniqueness of the Prandtl equations are proved
for initial data with multiple monotonicity regions,
by assuming that the initial data is tangentially real analytic on the complement of the monotonicity regions.
A typical example of the initial data in this category is:
\begin{align*}
\partial_2 u_{0,1}^P  < 0\quad {\rm for}~x_1<0\,, \quad \partial_2 u_{0,1}^P >0 \quad {\rm for}~x_1>0\,,
\qquad u_{0,1}^P ~{\rm is~real~analytic~in}~x_1~{\rm around}~x_1=0\,.
\end{align*}
That is, the monotonicity of the initial flow is lost around $x_1=0$, but instead,
the analyticity in the tangential variable is imposed there.
Because of this complementary distributions of two totally different structures and the nonlocal nature of the problem, the methods developed in the categories (I) and (II) are not enough in constructing local solutions in this class.
Indeed, there are several difficulties in this problem;
the norms for the analyticity class and the  monotonicity class are not compatible,
and moreover, nontrivial analytic functions cannot have a compact support,
which indicates the breakdown of standard localization arguments.
The key observation in \cite{KMVW14} is that one can in fact construct the analytic solution around $x_1=0$ in a decoupled manner without using the lateral boundary conditions. On the other hand, by introducing a suitable extension of the data in the tangential direction one can construct a monotone solution in the entire half plane.
Finally, by virtue of the finite propagation property in the tangential direction, the analytic solution and the monotonic solution actually coincide with each other on some strip regions, which implies the existence of solutions to the original Prandtl equations. The construction of \cite{KMVW14}  reveals in some sense a possibility of localizing the Prandtl equations in the tangential direction.
Moreover, the result of \cite{KMVW14} indicates that, even at the point of separation, the flow can be stable at least locally in time and space, if the flow is analytic in the tangential variable around the separation point.

\subsubsection{Ill-posedness results for the Prandtl equations}\label{subsubsec.illposed.Prandtl}

Although the local solvability of the unsteady Prandtl equations still remains open for general initial data in Sobolev spaces, several ill-posedness results have been reported in the literature.
This subsection is devoted to give an overview on recent progress in this direction.

\vspace{0.1cm}

\noindent
{\bf (I) Ill posedness of the Prandtl equations in Sobolev spaces.}
When given data are not monotonic, the unsteady Prandtl equations are known to be ill-posed in the sense of Hadamard.  The ill posedness is triggered by the instability of non-monotonic shear flows at high frequencies.
The first rigorous result for this instability is given in \cite{G-VD10}, where the linearization around a shear flow possessing a non-degenerate critical point is studied in details. To be precise, let $u^P=(u_1^P,u_2^P)$ be the solution to the Prandtl equations \eqref{Pra.sec3.4} for constant data $u^E=\underline{u}\in \mathbb{R}$ and $\partial_1p^E=0$ (that is, the Euler flow $u^0$ is a stationary shear flow $u^0=(u^0_1(x_2),0)$ and $\underline{u}=u^0_1 (0)$).
In this case, $u^P$ is also a shear flow $u^P(t,x) =(u_s(t,x_2),0)$, and $u_s$ is the solution to the heat equation:
\begin{equation}\label{heat.sec3.4}
  \left\{
\begin{aligned}
& \partial_t u_s   - \partial_2^2 u_s  = 0\,, \\
& u_s|_{t=0} = U_s \,, \qquad u_s|_{x_2=0} =0\,, \qquad  \lim_{x_2\rightarrow \infty} u_s = \underline{u}\,. \\
\end{aligned}\right.
\end{equation}
Here $U_s$ is a given initial shear profile satisfying the compatibility conditions.
Then, a natural question is whether or not one can construct a solution to the Prandtl equations around this shear flow. The key step to tackle this problem is to analyze the linearization around $u_s$:
\begin{equation}\label{LinearPra.sec3.4}
  \left\{
\begin{aligned}
& \partial_t v_1^P + u_s \partial_1 v_1^P + v_2^P \partial_2 u_s  - \partial_2^2 v_1^P  = 0 \,, \\
& v_2^P  = -\int_0^{x_2} \partial_{1} v_1^P  \dd y_2\,, \\
& v_1^P|_{t=s} = v^P_{0,1}\,, \qquad v_1^P|_{x_2=0} =0\,, \qquad  \lim_{x_2\rightarrow \infty} v_1^P = 0 \,. \\
\end{aligned}\right.
\end{equation}
Here, $t>s$, $x_1\in \mathbb{T}$, and $x_2\in \mathbb{R}_+$.
The system \eqref{LinearPra.sec3.4} is uniquely solvable at least locally in time, if the initial data $v_{0,1}^P$ is analytic in the $x_1$ variable, and then the evolution operator $T(t,s)$, $T(t,s) v_{0,1}^P:=v_1^P (t)$, is shown to be  locally well defined in the analytic functional framework; see \cite[Proposition 1]{G-VD10}.
With this observation one can define the operator norm of $T(t,s)$ from $H^{m_1}$ to $H^{m_2}$, where $H^{m_1}$ and $H^{m_2}$ are suitable Sobolev space in $\mathbb{T}\times \mathbb{R}_+$, and the exponents $m_1,m_2$ denote the order of the Sobolev regularity; see \cite{G-VD10} for the precise definition of $H^m$. In \cite{G-VD10}, the given initial data $U_s$ in \eqref{heat.sec3.4} is  assumed to have a no-degenerate critical point:
there is $a>0$ such that
\begin{align}\label{condition.critical}
U_s^{'}(a) =0\,, \qquad U_s^{''} (a) \ne 0\,.
\end{align}
Then the following ill posedness in the Sobolev class is given by \cite[Theorem 1]{G-VD10}.
\begin{thm}\label{thm.illposed.GVD} If \eqref{condition.critical} holds, then there exists $\sigma>0$ such that for all $\delta>0$,
\begin{align*}
\sup_{0\leq s\leq t\leq \delta} \| e^{-\sigma (t-s) \sqrt{|\partial_1|}} T (t,s) \|_{\mathcal{L}(H^{m},H^{m-\mu})} =\infty \qquad {\rm for~all}~~m\geq 0\,, \quad \mu\in [0,\frac12)\,.
\end{align*}
Moreover, there is a solution $u_s$ to \eqref{heat.sec3.4} and $\sigma>0$ such tat for all $\delta>0$,
\begin{align*}
\sup_{0\leq s\leq t\leq \delta} \| e^{-\sigma (t-s) \sqrt{|\partial_1|}} T (t,s) \|_{\mathcal{L}(H^{m_1},H^{m_2})} =\infty \qquad {\rm for~all}~~m_1,m_2\geq 0\,.
\end{align*}
\end{thm}
Th key step in the proof in \cite{G-VD10} is to construct an approximate solution to \eqref{LinearPra.sec3.4}, which grows in time with order $e^{\delta |n|^\frac12 (t-s)}$, $\delta>0$, for a tangential frequency $|n|\gg 1$. This growth rate $\mathcal{O}(|n|^\frac12)$ is responsible for the weight $e^{-\sigma (t-s)|\sqrt{|\partial_1|}}$ in the statement of Theorem \ref{thm.illposed.GVD}.
The approximate solution is constructed as a singular perturbation from an explicit solution to the inviscid linearized Prandtl equations (dropping the viscous term $\partial_2^2 v_1^P$ in \eqref{LinearPra.sec3.4}, and replacing $u_s$ by $U_s$),
for which the spectral problem has been studied in details in \cite{HoHu03}.
This instability mechanism, bifurcating from the inviscid solution and resulting in the growth rate $\mathcal{O}(|n|^\frac12)$, was first reported at a formal level in\cite{CoHoTu85}. Due to the nature of the singular perturbation, however, the rigorous justification of this mechanism requires a highly delicate asymptotic analysis, and it is successfully completed by \cite{G-VD10}.

The result of \cite{G-VD10} on the ill posedness for the linearized Prandtl equations is strengthened in \cite{GN11,G-VN12}, where it is shown that the solutions to the non-linear Prandtl equations cannot be Lipschitz continuous with respect to the initial data in Sobolev spaces.
Moreover, it is proved by \cite{GN11} that  in the Sobolev framework one cannot expect a natural estimate of the asymptotic boundary layer expansion for the Navier-Stokes flows, if the leading term in the boundary layer is a non-monotonic shear layer flow as in \cite{G-VD10}.

Recently the ill posedness of the three-dimensional Prandtl equations is studied in \cite{LWY16},
and it is revealed that there is a stronger instability mechanism in the three-dimensional case even at the linear level. In particular, it is shown in \cite{LWY16} that, in contrast to the two-dimensional case, the monotonicity condition on tangential velocity fields is not sufficient for the well posedness of the three-dimensional Prandtl equations.

\vspace{0.1cm}

\noindent
{\bf (II) Blowup solutions to the Prandtl equations.}
In the absence of the monotonicity, it is known that the solution to the Prandtl equations can blow up in  finite time. The formation
of such singularity was observed numerically in \cite{V-DoSh1980} for data corresponding to an impulsively started flow past a cylinder.  The existence of blowup solutions was also reported by \cite{HoHu03} through numerical and asymptotic analysis  for the inviscid Prandtl equations. The rigorous existence of finite-time blowup solutions is first given by \cite{EE97} (see also \cite{KVW16preprint}), which is stated as follows.
\begin{thm}\label{thm.blowup} Let $u^E=p^P=0$. Assume that the initial data $u_{0,1}^P$ is of the form $u_{0,1}^P (x_1,x_2) = - x_1 b_0 (x_1,x_2)$ for some smooth $b_0$, and that $a_0 (x_2) = - \partial_1 u_{0,1}^P (0,x_2)$ is non negative, smooth, and compactly supported. Assume in addition that
\begin{align}\label{negative.energy}
E (a_0) = \frac12 \| \partial_2 a_0 \|_{L^2(\mathbb{R}_+)}^2 - \frac14 \| a_0\|_{L^3 (\mathbb{R}_+)}^3 <0
\end{align}
holds. Then there exist no global smooth solutions to \eqref{Pra.sec3.4}.
\end{thm}

\begin{rem}{\rm In Theorem \ref{thm.blowup} the condition of the compact support of $a_0$ is not essential, and it can be replaced by a decay condition which ensures the boundedness of $\partial_1 u_1^P (t,0,x_2)$ in $L^1_{x_2} (\mathbb{R}_+)$ as long as the solution exists.
}
\end{rem}

A simple example of the initial data satisfying the conditions of Theorem \ref{thm.blowup} is
\begin{align}\label{example.blowup}
u_{0,1}^p (x_1,x_2) = -x_1 e^{-x_1^2} f (\frac{x_2}{R})\,,
\end{align}
where $f$ is a (non-trivial) non-negative smooth function with compact support, and $R>0$ is a sufficiently large number. Indeed, in this case the function $a_0$ is given by $a_0 (x_2) = - f(\frac{x_2}{R})$, and the quantity $E(a_0)$ is computed as
\begin{align*}
E (a_0) = \frac{1}{2R} \| \partial_2 f \|_{L^2 (\mathbb{R}_+)}^2 - \frac{R}{4} \| f\|_{L^3 (\mathbb{R}_+)}^3\,,
\end{align*}
which is negative when $R>0$ is large enough. Since the initial data defined by \eqref{example.blowup} is analytic in $x_1$, in virtue of Theorem \ref{thm.analyticity} there exists a unique solution to \eqref{Pra.sec3.4} at least for a short time. Theorem \ref{thm.blowup} shows that this local solution cannot be extended as a global solution, and the proof in \cite{EE97} implies that the blowup occurs for the quantity $\| \partial_1 u_1^P(t)\|_{L^\infty (\mathbb{R}^2_+)}$. The singularity formation observed in \cite{V-DoSh1980} is studied in details by \cite{GSS09}, in which numerical evidence is reported for the strong ill posedness of the Prandtl equations in the Sobolev space $H^1(\mathbb{R}_+)$.

\section{Conclusion}

The inviscid limit problem of the Navier-Stokes flows is one of the most fundamental and classical issues in fluid mechanics, particularly in understanding the flows at high Reynolds numbers.
Mathematically, the fundamental question here is whether or not the Navier-Stokes flows converge to the Euler flows in the zero-viscosity limit, by taking the effect of the boundary into account if necessary.
Even when there is no physical boundary, the analysis of the inviscid limit is a challenging problem if one works with singular flows such as vortex sheets or filaments (Section 2). The rigorous understanding of these structures, in terms of the analysis of the Navier-Stokes equations at high Reynolds numbers, is still out of reach except for the case when the distribution of the possible singularities can be well specified in advance due to some additional prescribed symmetry. In the presence of physical boundary, the verification of the inviscid limit is far from trivial in general even when given data have enough regularity, e.g., in higher order Sobolev spaces (Section 3). The main obstacle is the formation of the viscous boundary layer, the size and stability of which are crucially influenced by the type of the boundary conditions. If the boundary condition allows the flow to slip on the boundary or the boundary is non characteristic, then the effect of the boundary layer is moderate and the mathematical theory has been well developed by now in these categories (Sections 3.1, 3.3).  However, if the no-slip boundary condition is imposed, the size of the boundary layer is at least $\mathcal{O}(1)$ even at a formal level, and the underlying instability mechanism of the boundary layer leads to a serious difficulty in the analysis of the inviscid limit problem (Section 3.2).
The general result in this research area is Kato's criterion, which describes the condition for the convergence of the Navier-Stokes flows to the Euler flows in the energy space. This criterion can be confirmed under some symmetry conditions on both the fluid domain and the given data without assuming strong regularity of data such as analyticity.
Without symmetry, the known results verifying  the inviscid limit (the convergence in the energy space as well as the Prandtl boundary layer expansion) require, so far, the analyticity at least near the boundary, and the boundary is also assumed to be flat there. In understanding the formation of the boundary layer, mathematically the analysis of the Prandtl equations is a central issue (Section 3.4). However, the solvability of the Prandtl equations is available only for some restricted classes of given data including the monotonicity or the analyticity class (Section 3.4.1). In the general Sobolev framework, the strong instability and resulting ill posedness have been observed (Section 3.4.2). Moreover, without monotonicity large solutions to the Prandtl equations may blow up  in finite time, while the existence of smooth global solutions for small data, but in a wide fluid domain, remains open even under the monotonicity condition in general. Finally, even when the Prandtl equations are solvable, their solvability  does not imply the validity of the inviscid limit, and in fact, there is a significant discrepancy between these two problems.
This discrepancy actually indicates the limitation of Prandtl's approach to understand boundary layer separation in physical situations, although the latter has been often discussed within the framework of the Prandtl equations.
The mathematical understanding of the stability of the boundary layer and the analysis of the boundary layer separation needs a significant development of the present theory for the inviscid limit problem of the Navier-Stokes flows.

\section*{Cross references}

\begin{itemize}

\item Existence and stability of viscous vortices

\end{itemize}

\bibliographystyle{abbrv}


\begin{thebibliography}{100}

\bibitem{AbDa}
H.~Abidi and R.~Danchin.
\newblock Optimal bounds for the inviscid limit of {N}avier-{S}tokes equations.
\newblock {\em Asymptot. Anal.}, 38:35--46, 2004.

\bibitem{Aleks83}
S.~N. Alekseenko.
\newblock Solution of a degenerate linearized {N}avier-{S}tokes system with a
  homogeneous boundary condition.
\newblock In {\em Studies in integro-differential equations, {N}o. 16}, pages
  243--257. ``Ilim'', Frunze, 1983.

\bibitem{Aleks86}
S.~N. Alekseenko.
\newblock On vanishing viscosity in a linearized problem of the flow of an
  incompressible fluid.
\newblock In {\em Studies in integro-differential equations, {N}o.\ 19
  ({R}ussian)}, pages 288--304, 320. ``Ilim'', Frunze, 1986.

\bibitem{AWXY2015}
R.~Alexandre, Y.-G. Wang, C.-J. Xu, and T.~Yang.
\newblock Well-posedness of the prandtl equation in sobolev spaces.
\newblock {\em J. Amer. Math. Soc.}, 28(3):745--784, 2015.

\bibitem{Asano88}
K.~Asano.
\newblock A note on the abstract {C}auchy-{K}owalewski theorem.
\newblock {\em Proc. Japan Acad. Ser. A Math. Sci.}, 64(4):102--105, 1988.

\bibitem{Asano88'}
K.~Asano.
\newblock Zero-viscosity limit of the incompressible {N}avier-{S}tokes
  equations. {II}.
\newblock {\em Mathematical Analysis of Fluid and Plasma Dynamics,
  S${\bar{u}}$rikaisekikenky${\bar{u}}$sho K${\bar{o}}$ky${\bar{u}}$roku},
  656:105--128, 1988.

\bibitem{Bardos1972}
C.~Bardos.
\newblock Existence et unicit\'e de la solution de l'\'equation d'{E}uler en
  dimension deux.
\newblock {\em J. Math. Anal. Appl.}, 40:769--790, 1972.

\bibitem{BLNNT13}
C.~Bardos, M.~C. Lopes~Filho, D.~Niu, H.~J. Nussenzveig~Lopes, and E.~S. Titi.
\newblock Stability of two-dimensional viscous incompressible flows under
  three-dimensional perturbations and inviscid symmetry breaking.
\newblock {\em SIAM J. Math. Anal.}, 45(3):1871--1885, 2013.

\bibitem{BdVC10}
H.~Beir{\~a}o~da Veiga and F.~Crispo.
\newblock Sharp inviscid limit results under {N}avier type boundary conditions.
  {A}n {$L^p$} theory.
\newblock {\em J. Math. Fluid Mech.}, 12(3):397--411, 2010.

\bibitem{BdVC13}
H.~Beir{\~a}o~da Veiga and F.~Crispo.
\newblock Concerning the {$W^{k,p}$}-inviscid limit for 3-{D} flows under a
  slip boundary condition.
\newblock {\em J. Math. Fluid Mech.}, 13(1):117--135, 2011.

\bibitem{BdVC14}
H.~Beir{\~a}o~da Veiga and F.~Crispo.
\newblock The 3-{D} inviscid limit result under slip boundary conditions. {A}
  negative answer.
\newblock {\em J. Math. Fluid Mech.}, 14(1):55--59, 2012.

\bibitem{BdVC12}
H.~Beir{\~a}o~da Veiga and F.~Crispo.
\newblock A missed persistence property for the {E}uler equations and its
  effect on inviscid limits.
\newblock {\em Nonlinearity}, 25(6):1661--1669, 2012.

\bibitem{BNP04}
H.~Bellout, J.~Neustupa, and P.~Penel.
\newblock On the {N}avier-{S}tokes equation with boundary conditions based on
  vorticity.
\newblock {\em Math. Nachr.}, 269/270:59--72, 2004.

\bibitem{BS12}
L.~C. Berselli and S.~Spirito.
\newblock On the vanishing viscosity limit of 3{D} {N}avier-{S}tokes equations
  under slip boundary conditions in general domains.
\newblock {\em Comm. Math. Phys.}, 316(1):171--198, 2012.

\bibitem{BeCo}
A.~Bertozzi and P.~Constantin.
\newblock Global regularity for vortex patches.
\newblock {\em Comm. Math. Phys.}, 152(1):19--28, 1993.

\bibitem{BW02}
J.~L. Bona and J.~Wu.
\newblock The zero-viscosity limit of the 2{D} {N}avier-{S}tokes equations.
\newblock {\em Stud. Appl. Math.}, 109(4):265--278, 2002.

\bibitem{BruMar}
E.~Brunelli and C.~Marchioro.
\newblock Vanishing viscosity limit for a smoke ring with concentrated
  vorticity.
\newblock {\em J. Math. Fluid Mech.}, 13:421--428, 2011.

\bibitem{BFN10}
D.~Bucur, E.~Feireisl, and {\v{S}}.~Ne{\v{c}}asov{\'a}.
\newblock Boundary behavior of viscous fluids: influence of wall roughness and
  friction-driven boundary conditions.
\newblock {\em Arch. Ration. Mech. Anal.}, 197(1):117--138, 2010.

\bibitem{BILN12}
A.~V. Busuioc, D.~Iftimie, M.~C. Lopes~Filho, and H.~J. Nussenzveig~Lopes.
\newblock Incompressible {E}uler as a limit of complex fluid models with
  {N}avier boundary conditions.
\newblock {\em J. Differential Equations}, 252(1):624--640, 2012.

\bibitem{CaO1}
R.~E. Caflisch and O.~F. Orellana.
\newblock Long time existence for a slightly perturbed vortex sheet.
\newblock {\em Comm. Pure Appl. Math.}, 39(6):807--838, 1986.

\bibitem{CaO2}
R.~E. Caflisch and O.~F. Orellana.
\newblock Singular solutions and ill-posedness for the evolution of vortex
  sheets.
\newblock {\em SIAM J. Math. Anal.}, 20(2):293--307, 1989.

\bibitem{CaSa06}
R.~E. Caflisch and M.~Sammartino.
\newblock Vortex layers in the small viscosity limit.
\newblock In {\em WASCOM 2005-13th Conference on Waves and Stability in
  Continuous Media}, pages 59--70. World Scientific Publishing Company,
  Hackensack, 2006.

\bibitem{CLM13}
M.~Cannone, M.~C. Lombardo, and M.~Sammartino.
\newblock Well-posedness of {P}randtl equations with non-compatible data.
\newblock {\em Nonlinearity}, 26(12):3077--3100, 2013.

\bibitem{CLS14}
M.~Cannone, M.~C. Lombardo, and M.~Sammartino.
\newblock On the {P}randtl boundary layer equations in presence of corner
  singularities.
\newblock {\em Acta Appl. Math.}, 132:139--149, 2014.

\bibitem{ChaDu}
D.~Chae and P.~Dubovskii.
\newblock Functional and measure-valued solutions to the {E}uler equations for
  flows of incompressible fluids.
\newblock {\em Arch. Rational Mech. Anal.}, 129:385--396, 1995.

\bibitem{Che1}
J.-Y. Chemin.
\newblock Persistance de structures g{\'e}om{\'e}triques dans les fluides
  incompressibles bidimensionnels.
\newblock {\em Ann. Sci. Ecole Norm. Sup. (4)}, 26(4):517--542, 1993.

\bibitem{Che2}
J.-Y. Chemin.
\newblock A remark on the inviscid limit for two-dimensional incompressible
  fluids.
\newblock {\em Comm. Partial Differential Equations}, 21:1771--1779, 1996.

\bibitem{CMR98}
T.~Clopeau, A.~Mikeli{\'c}, and R.~Robert.
\newblock On the vanishing viscosity limit for the {$2{\rm D}$} incompressible
  {N}avier-{S}tokes equations with the friction type boundary conditions.
\newblock {\em Nonlinearity}, 11(6):1625--1636, 1998.

\bibitem{CKV15}
P.~Constantin, I.~Kukavica, and V.~Vicol.
\newblock On the inviscid limit of the {N}avier-{S}tokes equations.
\newblock {\em Proc. Amer. Math. Soc.}, 143(7):3075--3090, 2015.

\bibitem{CoWu1}
P.~Constantin and J.~Wu.
\newblock Inviscid limit for vortex patches.
\newblock {\em Nonlinearity}, 8:735--742, 1995.

\bibitem{CoWu2}
P.~Constantin and J.~Wu.
\newblock The inviscid limit for non-smooth vorticity.
\newblock {\em Indiana Univ. Math. J.}, 45:67--81, 1996.

\bibitem{CoHoTu85}
S.~J. Cowley, L.~M. Hocking, and O.~R. Tutty.
\newblock The stability of solutions of the classical unsteady boundary-layer
  equation.
\newblock {\em Phys. Fluids}, 28:441--443, 1985.

\bibitem{Dan1}
R.~Danchin.
\newblock Poches de tourbillon visqueuses.
\newblock {\em J. Math. Pures Appl.}, 76:609--647, 1997.

\bibitem{Dan2}
R.~Danchin.
\newblock Persistance de structures g{\'e}om{\'e}triques et limite non
  visqueuse pour les fluides incompressibles en dimension quelconque.
\newblock {\em Bull. Soc. Math. France}, 127:179--227, 1999.

\bibitem{De}
N.~Depauw.
\newblock Poche de tourbillon pour {E}uler $2${D} dans un ouvert {\`a} bord.
\newblock {\em J. Math. Pures Appl. (9)}, 78(3):313--351, 1999.

\bibitem{DuRo}
J.~Duchon and R.~Robert.
\newblock Global vortex sheet solutions of {E}uler equations in the plane.
\newblock {\em J. Differential Equations}, 73(2):215--224, 1988.

\bibitem{Du}
A.~Dutrifoy.
\newblock On $3$-{D} vortex patches in bounded domains.
\newblock {\em Comm. Partial Differential Equations}, 28(7-8):1237--1263, 2003.

\bibitem{EE97}
W.~E and B.~Engquist.
\newblock Blowup of solutions of the unsteady {P}randtl's equation.
\newblock {\em Comm. Pure Appl. Math.}, 50(12):1287--1293, 1997.

\bibitem{EbMa}
D.~Ebin and J.~Marsden.
\newblock Groups of diffeomorphisms and the notion of an incompressible fluid.
\newblock {\em Ann. Math.}, 92:102--163, 1970.

\bibitem{FTZpreprint}
M.~Fei, T.~Tao, and Z.~Zhang.
\newblock On the zero-viscosity limit of the {N}avier-{S}tokes equations in the
  half-space.
\newblock {\em ArXiv e-prints}, Sept. 2016.

\bibitem{FeSve}
H.~Feng and V.~{\v S}ver\'ak.
\newblock On the {C}auchy problem for axi-symmetric vortex rings.
\newblock {\em Arch. Rational Mech. Anal.}, 215:89--123, 2015.

\bibitem{GaGa}
I.~Gallagher and T.~Gallay.
\newblock Uniqueness for the two-dimensional {N}avier-{S}tokes equation with a
  measure as initial vorticity.
\newblock {\em Math. Ann.}, 332:287--327, 2005.

\bibitem{Gal}
T.~Gallay.
\newblock Interaction of vortices in weakly viscous planar flows.
\newblock {\em Arch. Rational Mech. Anal.}, 200:445--490, 2011.

\bibitem{GaMa}
T.~Gallay and Y.~Maekawa.
\newblock Existence and stability of viscous vortices.
\newblock In {\em Handbook of Mathematical Analysis of Mechanics in Viscous
  Fluids}. Springer.

\bibitem{GaWa}
T.~Gallay and C.~E. Wayne.
\newblock Global stability of vortex solutions of the two-dimensional
  {N}avier-{S}tokes equation.
\newblock {\em Comm. Math. Phys.}, 255:97--129, 2005.

\bibitem{GaSa}
P.~Gamblin and X.~Saint~Raymond.
\newblock On three-dimensional vortex patches.
\newblock {\em Bull. Soc. Math. France}, 123(3):375--424, 1995.

\bibitem{GSS09}
F.~Gargano, M.~Sammartino, and V.~Sciacca.
\newblock Singularity formation for {P}randtl's equations.
\newblock {\em Phys. D}, 238(19):1975--1991, 2009.

\bibitem{G-VD10}
D.~G{\'e}rard-Varet and E.~Dormy.
\newblock On the ill-posedness of the {P}randtl equation.
\newblock {\em J. Amer. Math. Soc.}, 23(2):591--609, 2010.

\bibitem{GMMpreprint}
D.~G{\'e}rard-Varet, Y.~Maekawa, and N.~Masmoudi.
\newblock Gevrey stability of {P}randtl expansions for $2${D} {N}avier-{S}tokes
  flows.
\newblock {\em ArXiv e-prints}, July 2016.

\bibitem{D-VM15}
D.~G{\'e}rard-Varet and N.~Masmoudi.
\newblock Well-posedness for the {P}randtl system without analyticity or
  monotonicity.
\newblock {\em Ann. Scient. {\'E}c. Norm. Sup.}, 48(4):1273--1325, 2015.

\bibitem{G-VN12}
D.~G{\'e}rard-Varet and T.~Nguyen.
\newblock Remarks on the ill-posedness of the {P}randtl equation.
\newblock {\em Asymptot. Anal.}, 77(1-2):71--88, 2012.

\bibitem{Gie13}
G.-M. Gie.
\newblock Asymptotic expansion of the {S}tokes solutions at small viscosity:
  the case of non-compatible initial data.
\newblock {\em Commun. Math. Sci.}, 12(2):383--400, 2014.

\bibitem{GJ13}
G.-M. Gie and C.-Y. Jung.
\newblock Vorticity layers of the 2{D} {N}avier-{S}tokes equations with a slip
  type boundary condition.
\newblock {\em Asymptot. Anal.}, 84(1-2):17--33, 2013.

\bibitem{GK12}
G.-M. Gie and J.~P. Kelliher.
\newblock Boundary layer analysis of the {N}avier-{S}tokes equations with
  generalized {N}avier boundary conditions.
\newblock {\em J. Differential Equations}, 253(6):1862--1892, 2012.

\bibitem{GKLMN15}
G.-M. Gie, J.~P. Kelliher, M.~C. Lopes, A.~L. Mazzucato, and H.~J.
  Nussenzveig~Lopes.
\newblock The vanishing viscosity limit for some symmetric flows.
\newblock In final preparation.

\bibitem{GiMi}
Y.~Giga and T.~Miyakawa.
\newblock {N}avier-{S}tokes flow in $\mathbb{R}^3$ with measures as initial
  vorticity and {M}orrey spaces.
\newblock {\em Comm. Partial Differential Equations}, 14(5):577--618, 1989.

\bibitem{GiMiO}
Y.~Giga, T.~Miyakawa, and H.~Osada.
\newblock Two-dimensional {N}avier-{S}tokes flow with measures as initial
  vorticity.
\newblock {\em Arch. Rational Mech. Anal.}, 104:223--250, 1988.

\bibitem{Go}
K.~K. Golovkin.
\newblock Vanishing viscosity in the cauchy problem for equations of
  hydrodynamics.
\newblock {\em Trudy Mat. Inst. Steklov.}, 92:31--49, 1966.

\bibitem{Grenier00}
E.~Grenier.
\newblock On the nonlinear instability of {E}uler and {P}randtl equations.
\newblock {\em Comm. Pure Appl. Math.}, 53(9):1067--1091, 2000.

\bibitem{GGNpreprint'}
E.~Grenier, Y.~Guo, and T.~Nguyen.
\newblock Spectral instability of characteristic boundary layer flows.
\newblock {\em ArXiv e-prints}, June 2014.

\bibitem{GN11}
Y.~Guo and T.~Nguyen.
\newblock A note on {P}randtl boundary layers.
\newblock {\em Comm. Pure Appl. Math.}, 64(10):1416--1438, 2011.

\bibitem{HMNW11}
D.~Han, A.~L. Mazzucato, D.~Niu, and X.~Wang.
\newblock Boundary layer for a class of nonlinear pipe flow.
\newblock {\em J. Differential Equations}, 252(12):6387--6413, 2012.

\bibitem{HW14}
D.~Han and X.~Wang.
\newblock Initial-boundary layer associated with the nonlinear
  {D}arcy-{B}rinkman system.
\newblock {\em J. Differential Equations}, 256(2):609--639, 2014.

\bibitem{Hmi1}
T.~Hmidi.
\newblock R{\'e}gularit{\'e} h{\''o}ld{\'e}rienne des poches de tourbillon
  visqueuses.
\newblock {\em J. Math. Pures Appl.}, 84:1455--1495, 2005.

\bibitem{Hmi2}
T.~Hmidi.
\newblock Poches de tourbillon singuli{\`e}res dans un fluide faiblement
  visqueux.
\newblock {\em Rev. Mat. Iberoamericana}, 22:489--543, 2006.

\bibitem{HoHu03}
L.~Hong and J.~Hunter.
\newblock Singularity formation and instability in the unsteady inviscid and
  viscous {P}randtl equations.
\newblock {\em Comm. Math. Sci.}, 1:293--316, 2003.

\bibitem{Hua1}
C.~Huang.
\newblock Remarks on regularity of non-constant vortex patches.
\newblock {\em Commun. Appl. Anal.}, 3(4):449--459, 1999.

\bibitem{Hua2}
C.~Huang.
\newblock Singular integral system approach to regularity of $3$d vortex
  patches.
\newblock {\em Indiana Univ. Math. J.}, 50(1):509--552, 2001.

\bibitem{ILN09}
D.~Iftimie, M.~C. Lopes~Filho, and H.~J. Nussenzveig~Lopes.
\newblock Incompressible flow around a small obstacle and the vanishing
  viscosity limit.
\newblock {\em Comm. Math. Phys.}, 287(1):99--115, 2009.

\bibitem{IP06}
D.~Iftimie and G.~Planas.
\newblock Inviscid limits for the {N}avier-{S}tokes equations with {N}avier
  friction boundary conditions.
\newblock {\em Nonlinearity}, 19(4):899--918, 2006.

\bibitem{IS11}
D.~Iftimie and F.~Sueur.
\newblock Viscous boundary layers for the {N}avier-{S}tokes equations with the
  {N}avier slip conditions.
\newblock {\em Arch. Ration. Mech. Anal.}, 199(1):145--175, 2011.

\bibitem{IV2016}
M.~Ignatova and V.~Vicol.
\newblock Almost global existence for the {P}randtl boundary layer equations.
\newblock {\em Arch. Rational Mech. Anal.}, 220:809--848, 2016.

\bibitem{JM01}
W.~J{\"a}ger and A.~Mikeli{\'c}.
\newblock On the roughness-induced effective boundary conditions for an
  incompressible viscous flow.
\newblock {\em J. Differential Equations}, 170(1):96--122, 2001.

\bibitem{Kato1972}
T.~Kato.
\newblock Nonstationary flows of viscous and ideal fluids in {${\bf R}\sp{3}$}.
\newblock {\em J. Functional Analysis}, 9:296--305, 1972.

\bibitem{K84}
T.~Kato.
\newblock Remarks on zero viscosity limit for nonstationary {N}avier-{S}tokes
  flows with boundary.
\newblock In {\em Seminar on nonlinear partial differential equations
  ({B}erkeley, {C}alif., 1983)}, volume~2 of {\em Math. Sci. Res. Inst. Publ.},
  pages 85--98. Springer, New York, 1984.

\bibitem{Ka3}
T.~Kato.
\newblock The {N}avier-{S}tokes equation for an incompressible fluid in
  $\mathbb{R}^2$ with a measure as the initial vorticity.
\newblock {\em Differ. Integral Equ.}, 7:949--966, 1994.

\bibitem{K06}
J.~P. Kelliher.
\newblock Navier-{S}tokes equations with {N}avier boundary conditions for a
  bounded domain in the plane.
\newblock {\em SIAM J. Math. Anal.}, 38(1):210--232 (electronic), 2006.

\bibitem{K07}
J.~P. Kelliher.
\newblock On {K}ato's conditions for vanishing viscosity.
\newblock {\em Indiana Univ. Math. J.}, 56(4):1711--1721, 2007.

\bibitem{K08}
J.~P. Kelliher.
\newblock Vanishing viscosity and the accumulation of vorticity on the
  boundary.
\newblock {\em Commun. Math. Sci.}, 6(4):869--880, 2008.

\bibitem{K14}
J.~P. {Kelliher}.
\newblock {Observations on the vanishing viscosity limit}.
\newblock {\em ArXiv e-prints}, Sept. 2014.

\bibitem{KLN09}
J.~P. Kelliher, M.~C.~L. Filho, and H.~J.~N. Lopes.
\newblock Vanishing viscosity limit for an expanding domain in space.
\newblock {\em Ann. Inst. H. Poincar\'e Anal. Non Lin\'eaire},
  26(6):2521--2537, 2009.

\bibitem{KTW11}
J.~P. Kelliher, R.~Temam, and X.~Wang.
\newblock Boundary layer associated with the {D}arcy-{B}rinkman-{B}oussinesq
  model for convection in porous media.
\newblock {\em Phys. D}, 240(7):619--628, 2011.

\bibitem{KMVW14}
I.~Kukavica, N.~Masmoudi, V.~Vicol, and T.~K. Wong.
\newblock On the local well-posedness of the {P}randtl and hydrostatic {E}uler
  equations with multiple monotonicity regions.
\newblock {\em SIAM J. Math. Anal.}, 46(6):3865--3890, 2014.

\bibitem{KV13}
I.~Kukavica and V.~Vicol.
\newblock On the local existence of analytic solutions to the {P}randtl
  boundary layer equations.
\newblock {\em Commun. Math. Sci.}, 11(1):269--292, 2013.

\bibitem{KVW16preprint}
I.~{Kukavica}, V.~{Vicol}, and F.~{Wang}.
\newblock {The van Dommelen and Shen singularity in the Prandtl equations}.
\newblock {\em ArXiv e-prints}, Dec. 2015.

\bibitem{LCM15}
C.~Lacave and A.~L. Mazzucato.
\newblock The vanishing viscosity limit in the presence of a porous medium.
\newblock {\em Math. Ann.}, 2015.

\bibitem{LBS07}
E.~Lauga, M.~Brenner, and H.~Stone.
\newblock Microfluidics: The no-slip boundary condition.
\newblock In C.~Tropea, A.~Yarin, and J.~F. Foss, editors, {\em Springer
  Handbook of Experimental Fluid Mechanics}, pages 1219--1240. Springer, 2007.

\bibitem{Le}
G.~Lebeau.
\newblock R{\'e}gularit{\'e} du probl{\`e}me de {K}elvin-{H}elmholtz pour
  l'{\'e}quation d'euler $2${D}.
\newblock In {\em S{\'e}minaire: Equations aux D{\'e}riv{\'e}es Partielles,
  2000-2001, Exp. No. II, S{\'e}minaire: {\'E}quations aux D{\'e}riv{\'e}es
  Partielles}, page 12 pp. {\'E}cole Polytechnique, Palaiseau, France, 2001.

\bibitem{Lions69}
J.-L. Lions.
\newblock {\em Quelques m\'ethodes de r\'esolution des probl\`emes aux limites
  non lin\'eaires}.
\newblock Dunod; Gauthier-Villars, Paris, 1969.

\bibitem{Lions73}
J.-L. Lions.
\newblock {\em Perturbations singuli\`eres dans les probl\`emes aux limites et
  en contr\^ole optimal}.
\newblock Lecture Notes in Mathematics, Vol. 323. Springer-Verlag, Berlin-New
  York, 1973.

\bibitem{L1996}
P.-L. Lions.
\newblock {\em Mathematical topics in fluid mechanics. {V}ol. 1}, volume~3 of
  {\em Oxford Lecture Series in Mathematics and its Applications}.
\newblock The Clarendon Press Oxford University Press, New York, 1996.

\bibitem{LiWa16}
C.-J. Liu and Y.-G. Wang.
\newblock Derivation of {P}randtl boundary layer equations for the
  incompressible {N}avier-{S}tokes equations in a curved domain.
\newblock {\em Appl. Math. Lett.}, 34:81--85, 2014.

\bibitem{LWYPreprint}
C.~J. Liu, Y.-G. Wang, and T.~Yang.
\newblock {A well-posedness theory for the Prandtl equations in three space
  variables}.
\newblock {\em ArXiv e-prints}, May 2014.

\bibitem{LWY16}
C.-J. Liu, Y.-G. Wang, and T.~Yang.
\newblock On the ill-posedness of the {P}randtl equations in three-dimensional
  space.
\newblock {\em Arch. Rational Mech. Anal.}, 220:83--108, 2016.

\bibitem{LCS03}
M.~C. Lombardo, M.~Cannone, and M.~Sammartino.
\newblock Well-posedness of the boundary layer equations.
\newblock {\em SIAM J. Math. Anal.}, 35(4):987--1004 (electronic), 2003.

\bibitem{LSOseen01}
M.~C. Lombardo and M.~Sammartino.
\newblock Zero viscosity limit of the {O}seen equations in a channel.
\newblock {\em SIAM J. Math. Anal.}, 33(2):390--410 (electronic), 2001.

\bibitem{LMN08}
M.~C. Lopes~Filho, A.~L. Mazzucato, and H.~J. Nussenzveig~Lopes.
\newblock Vanishing viscosity limit for incompressible flow inside a rotating
  circle.
\newblock {\em Phys. D}, 237(10-12):1324--1333, 2008.

\bibitem{LMNT08}
M.~C. Lopes~Filho, A.~L. Mazzucato, H.~J. Nussenzveig~Lopes, and M.~Taylor.
\newblock Vanishing viscosity limits and boundary layers for circularly
  symmetric 2{D} flows.
\newblock {\em Bull. Braz. Math. Soc. (N.S.)}, 39(4):471--513, 2008.

\bibitem{LNP05}
M.~C. Lopes~Filho, H.~J. Nussenzveig~Lopes, and G.~Planas.
\newblock On the inviscid limit for two-dimensional incompressible flow with
  {N}avier friction condition.
\newblock {\em SIAM J. Math. Anal.}, 36(4):1130--1141 (electronic), 2005.

\bibitem{LoNuScho}
M.~C. Lopes~Filho, H.~J. Nussenzveig~Lopes, and S.~Schochet.
\newblock A criterion for the equivalence of the {B}irkhoff-{R}ott and {E}uler
  descriptions of vortex sheet evolution.
\newblock {\em Trans. Amer. Math. Soc.}, 359(9):4125--4142, 2007.

\bibitem{LMNTZ15}
M.~C. Lopes~Filho, H.~J. Nussenzveig~Lopes, E.~S. Titi, and A.~Zang.
\newblock Approximation of 2{D} {E}uler equations by the second-grade fluid
  equations with {D}irichlet boundary conditions.
\newblock {\em J. Math. Fluid Mech.}, 17(2):327--340, 2015.

\bibitem{M14}
Y.~Maekawa.
\newblock On the inviscid limit problem of the vorticity equations for viscous
  incompressible flows in the half-plane.
\newblock {\em Comm. Pure Appl. Math.}, 67(7):1045--1128, 2014.

\bibitem{Maj}
A.~Majda.
\newblock Vorticity and the mathematical theory of incompressible fluid flow.
\newblock {\em Comm. Pure Appl. Math.}, 39:S187--S220, 1986.

\bibitem{MajdaBertozzi}
A.~J. Majda and A.~L. Bertozzi.
\newblock {\em Vorticity and incompressible flow}, volume~27 of {\em Cambridge
  Texts in Applied Mathematics}.
\newblock Cambridge University Press, Cambridge, 2002.

\bibitem{Mar1}
C.~Marchioro.
\newblock On the vanishing viscosity limit for two-dimensional
  {N}avier-{S}tokes equations with singular initial data.
\newblock {\em Math. Methods Appl. Sci.}, 12:463--470, 1990.

\bibitem{Mar2}
C.~Marchioro.
\newblock On the inviscid limit for a fluid with a concentrated vorticity.
\newblock {\em Comm. Math. Phys.}, 196:53--65, 1998.

\bibitem{Mar3}
C.~Marchioro.
\newblock Vanishing viscosity limit for an incompressible fluid with
  concentrated vorticity.
\newblock {\em J. Math. Phys.}, 48:065302, 2007.

\bibitem{MarPu}
C.~Marchioro and M.~Pulvirenti.
\newblock Vortices and localization in {E}uler flows.
\newblock {\em Comm. Math. Phys.}, 154:49--61, 1993.

\bibitem{Mas2007}
N.~Masmoudi.
\newblock Examples of singular limits in hydrodynamics.
\newblock In {\em Handbook of differential equations: evolutionary equations.
  Vol III}, pages 195--275. Elsevier/North-Holland, Amsterdam, 2007.

\bibitem{Masmoudi2007}
N.~Masmoudi.
\newblock Remarks about the inviscid limit of the {N}avier-{S}tokes system.
\newblock {\em Comm. Math. Phys.}, 270(3):777--788, 2007.

\bibitem{MR12}
N.~Masmoudi and F.~Rousset.
\newblock Uniform regularity for the {N}avier-{S}tokes equation with {N}avier
  boundary condition.
\newblock {\em Arch. Ration. Mech. Anal.}, 203(2):529--575, 2012.

\bibitem{MaWo2015}
N.~Masmoudi and T.~Wong.
\newblock Local-in-time existence and uniqueness of solutions to the prandtl
  equations by energy methods.
\newblock {\em Comm. Pure Appl. Math.}, 68:1683--1741, 2015.

\bibitem{Matsui1994}
S.~Matsui.
\newblock Example of zero viscosity limit for two-dimensional nonstationary
  {N}avier-{S}tokes flows with boundary.
\newblock {\em Japan J. Indust. Appl. Math.}, 11(1):155--170, 1994.

\bibitem{MaShi1984}
S.~Matsui and T.~Shirota.
\newblock On separation points of solutions to {P}randtl boundary layer
  problem.
\newblock {\em Hokkaido Math. J.}, 13:92--108, 1984.

\bibitem{MNW10}
A.~Mazzucato, D.~Niu, and X.~Wang.
\newblock Boundary layer associated with a class of 3{D} nonlinear plane
  parallel channel flows.
\newblock {\em Indiana Univ. Math. J.}, 60(4):1113--1136, 2011.

\bibitem{MT11}
A.~Mazzucato and M.~Taylor.
\newblock Vanishing viscosity limits for a class of circular pipe flows.
\newblock {\em Comm. Partial Differential Equations}, 36(2):328--361, 2011.

\bibitem{MT08}
A.~L. Mazzucato and M.~E. Taylor.
\newblock Vanishing viscosity plane parallel channel flow and related singular
  perturbation problems.
\newblock {\em Analysis $\&$ PDE}, 1(1):35--93, 2008.

\bibitem{Mc}
F.~J. McGrath.
\newblock Nonstationary plane flow of viscous and ideal fluids.
\newblock {\em Arch. Rational Mech. Anal.}, 27:329--348, 1968.

\bibitem{MNN13}
A.~Mikeli{\'c}, {\v{S}}.~Ne{\v{c}}asov{\'a}, and M.~Neuss-Radu.
\newblock Effective slip law for general viscous flows over an oscillating
  surface.
\newblock {\em Math. Methods Appl. Sci.}, 36(15):2086--2100, 2013.

\bibitem{MikelicPaoli}
A.~Mikeli{\'c} and L.~Paoli.
\newblock Homogenization of the inviscid incompressible fluid flow through a
  {$2$}{D} porous medium.
\newblock {\em Proc. Amer. Math. Soc.}, 127(7):2019--2028, 1999.

\bibitem{Oleinik63}
O.~A. Ole{\u\i}nik.
\newblock On the system of equations of the boundary layer theory.
\newblock {\em Zh. bychisl. matem. i matem. fiz.}, 3(3):489--507, 1963.

\bibitem{Oleinik66}
O.~A. Ole{\u\i}nik.
\newblock On the mathematical theory of boundary layer for an unsteady flow of
  incompressible fluid.
\newblock {\em J. Appl. Math. Mech.}, 30:951--974 (1967), 1966.

\bibitem{Oleinik66Stability}
O.~A. Ole{\u\i}nik.
\newblock On the stability of solutions of the system of boundary layer
  equations for a nonstationary flow of an incompressible fluid.
\newblock {\em PMM}, 30(3):417--423, 1966.

\bibitem{OleinikSamokhin99}
O.~A. Ole{\u\i}nik and V.~N. Samokhin.
\newblock {\em Mathematical models in boundary layer theory}, volume~15 of {\em
  Applied Mathematics and Mathematical Computation}.
\newblock Chapman \& Hall/CRC, Boca Raton, FL, 1999.

\bibitem{Pa14}
M.~Paddick.
\newblock Stability and instability of {N}avier boundary layers.
\newblock {\em Differ. Integral Equ.}, 27:893--930, 2014.

\bibitem{Pra1904}
L.~Prandtl.
\newblock {${\rm \ddot{U}}$}ber {F}l{${\rm \ddot{u}}$}ssigkeitsbewegung bei
  sehr kleiner {R}eibung.
\newblock {\em Verh. III Intern. Math. Kongr. Heidelberg}, pages 485--491,
  1904.

\bibitem{Ru06}
W.~M. Rusin.
\newblock On the inviscid limit for the solutions of two-dimensional
  incompressible {N}avier-{S}tokes equations with slip-type boundary
  conditions.
\newblock {\em Nonlinearity}, 19(6):1349--1363, 2006.

\bibitem{SC98}
M.~Sammartino and R.~E. Caflisch.
\newblock Zero viscosity limit for analytic solutions of the {N}avier-{S}tokes
  equation on a half-space. {I}, {II}.
\newblock {\em Comm. Math. Phys.}, 192(2):433--491, 1998.

\bibitem{Schlichting}
H.~Schlichting and K.~Gersten.
\newblock {\em Boundary-layer theory}.
\newblock Springer-Verlag, Berlin, enlarged edition, 2000.
\newblock With contributions by Egon Krause and Herbert Oertel, Jr., Translated
  from the ninth German edition by Katherine Mayes.

\bibitem{S12}
F.~Sueur.
\newblock A {K}ato type theorem for the inviscid limit of the {N}avier-{S}tokes
  equations with a moving rigid body.
\newblock {\em Comm. Math. Phys.}, 316(3):783--808, 2012.

\bibitem{S14}
F.~Sueur.
\newblock On the inviscid limit for the compressible {N}avier-{S}tokes system
  in an impermeable bounded domain.
\newblock {\em J. Math. Fluid Mech.}, 16(1):163--178, 2014.

\bibitem{Sue}
F.~Sueur.
\newblock Viscous profiles of vortex patches.
\newblock {\em J. Inst. Math. Jussieu}, 14(1):1--68, 2015.

\bibitem{SuSuBarFri}
C.~Sulem, P.-L. Sulem, C.~Bardos, and U.~Frisch.
\newblock Finite time analyticity for the two- and three-dimensional
  {K}elvin-{H}elmholtz instability.
\newblock {\em Comm. Math. Phys.}, 80(4):485--516, 1981.

\bibitem{Swann1971}
H.~S.~G. Swann.
\newblock The convergence with vanishing viscosity of nonstationary
  {N}avier-{S}tokes flow to ideal flow in {$R\sb{3}$}.
\newblock {\em Trans. Amer. Math. Soc.}, 157:373--397, 1971.

\bibitem{TWOseen96}
R.~Temam and X.~Wang.
\newblock Asymptotic analysis of {O}seen type equations in a channel at small
  viscosity.
\newblock {\em Indiana Univ. Math. J.}, 45(3):863--916, 1996.

\bibitem{TWOseen97}
R.~Temam and X.~Wang.
\newblock Boundary layers for {O}seen's type equation in space dimension three.
\newblock {\em Russian J. Math. Phys.}, 5(2):227--246 (1998), 1997.

\bibitem{TW1998}
R.~Temam and X.~Wang.
\newblock On the behavior of the solutions of the {N}avier-{S}tokes equations
  at vanishing viscosity.
\newblock {\em Ann. Scuola Norm. Sup. Pisa Cl. Sci. (4)}, 25(3-4):807--828
  (1998), 1997.
\newblock Dedicated to Ennio De Giorgi.

\bibitem{TWsuction}
R.~Temam and X.~Wang.
\newblock Boundary layers in channel flow with injection and suction.
\newblock {\em Appl. Math. Lett.}, 14(1):87--91, 2001.

\bibitem{TW02}
R.~Temam and X.~Wang.
\newblock Boundary layers associated with incompressible {N}avier-{S}tokes
  equations: the noncharacteristic boundary case.
\newblock {\em J. Differential Equations}, 179(2):647--686, 2002.

\bibitem{V-DoSh1980}
L.~L. Van~Dommelen and S.~F. Shen.
\newblock The spontaneous generation of the singularity in a separating laminar
  boundary layer.
\newblock {\em J. Comput. Phys.}, 38:125--140, 1980.

\bibitem{VanDyke82}
M.~Van~Dyke.
\newblock {\em An album of fluid motion}.
\newblock Parabolic Press, 1982.

\bibitem{Vi}
M.~Vishik.
\newblock Incompressible flows of an ideal fluid with vorticity in borderline
  spaces of {B}esov type.
\newblock {\em Ann. Sci. Ecole Norm. Sup.}, 32:769--812, 1999.

\bibitem{VL57}
M.~I. Vishik and L.~A. Lyusternik.
\newblock Regular degeneration and boundary layer for linear differential
  equations with small parameter.
\newblock {\em Uspekki Mat. Nauk}, 12:3--122, 1957.

\bibitem{W2001}
X.~Wang.
\newblock A {K}ato type theorem on zero viscosity limit of {N}avier-{S}tokes
  flows.
\newblock {\em Indiana Univ. Math. J.}, 50(Special Issue):223--241, 2001.
\newblock Dedicated to Professors Ciprian Foias and Roger Temam (Bloomington,
  IN, 2000).

\bibitem{WX05}
Y.-G. Wang and Z.~Xin.
\newblock Zero-viscosity limit of the linearized compressible {N}avier-{S}tokes
  equations with highly oscillatory forces in the half-plane.
\newblock {\em SIAM J. Math. Anal.}, 37(4):1256--1298, 2005.

\bibitem{Wu}
S.~Wu.
\newblock Mathematical analysis of vortex sheets.
\newblock {\em Comm. Pure Appl. Math.}, 59(8):1065--1206, 2006.

\bibitem{XiaoXin07}
Y.~Xiao and Z.~Xin.
\newblock On the vanishing viscosity limit for the 3{D} {N}avier-{S}tokes
  equations with a slip boundary condition.
\newblock {\em Comm. Pure Appl. Math.}, 60(7):1027--1055, 2007.

\bibitem{XXW09}
Y.~Xiao, Z.~Xin, and J.~Wu.
\newblock Vanishing viscosity limit for the 3{D} magnetohydrodynamic system
  with a slip boundary condition.
\newblock {\em J. Funct. Anal.}, 257(11):3375--3394, 2009.

\bibitem{XL11}
X.~Xie and C.~Li.
\newblock Vanishing viscosity limit for viscous magnetohydrodynamic equations
  with a slip boundary condition.
\newblock {\em Appl. Math. Sci. (Ruse)}, 5(41-44):1999--2011, 2011.

\bibitem{XY99}
Z.~Xin and T.~Yanagisawa.
\newblock Zero-viscosity limit of the linearized {N}avier-{S}tokes equations
  for a compressible viscous fluid in the half-plane.
\newblock {\em Comm. Pure Appl. Math.}, 52(4):479--541, 1999.

\bibitem{XZ04}
Z.~Xin and L.~Zhang.
\newblock On the global existence of solutions to the {P}randtl's system.
\newblock {\em Adv. Math.}, 181:88--133, 2004.

\bibitem{Y1963}
V.~I. Yudovich.
\newblock Non-stationary flows of an ideal incompressible fluid.
\newblock {\em \u Z. Vy\v cisl. Mat. i Mat. Fiz.}, 3:1032--1066 (Russian),
  1963.

\bibitem{Yu2}
V.~I. Yudovich.
\newblock Uniqueness theorem for the basic nonstationary problem in the
  dynamics of an ideal incompressible fluid.
\newblock {\em Math. Res. Lett.}, 2:27--38, 1995.

\bibitem{ZZ2016}
P.~Zhang and Z.~Zhang.
\newblock Long time well-posedness of {P}randtl system with small and analytic
  initial data.
\newblock {\em Journal of Functional Analysis}, 270:2591--2615, 2016.

\end{thebibliography}

\end{document}